\documentclass{article}

\PassOptionsToPackage{compress}{natbib}

     \usepackage[preprint]{neurips_2021}

\usepackage[utf8]{inputenc} 
\usepackage[T1]{fontenc}    
\usepackage{hyperref}       
\usepackage{url}            
\usepackage{booktabs}       
\usepackage{amsfonts}       
\usepackage{nicefrac}       
\usepackage{microtype}      
\usepackage{xcolor}         

\title{Stochastic Learning Rate Optimization in the Stochastic Approximation and Online Learning Settings}

\author{Theodoros~Mamalis \\
	Department of Electrical and Computer Engineering\\
	University of Illinois at Urbana-Champaign\\
	306 N Wright St, Urbana, IL 61801, USA \\
	\texttt{mamalis2@illinois.edu} \\
	\And
	Du{\v s}an Stipanovi{\' c} \\
	Coordinated Science Laboratory\\
	University of Illinois at Urbana-Champaign\\
	1308 W Main St, Urbana, IL 61801, USA \\
	\texttt{dusan@illinois.edu} \\
	\AND
	Petros Voulgaris \\
	Department of Mechanical Engineering\\
	University of Nevada\\
	Reno, NV 89557, USA \\
	\texttt{pvoulgaris@unr.edu} \\
}

\usepackage{url}

\usepackage{caption}
\usepackage{subcaption}
\usepackage{graphicx}

\usepackage{amsmath}
\usepackage{amssymb} 
\usepackage{amsthm}

\usepackage{balance} 

\newcommand{\E}[2][]{{\mathbb{E}_{#1}[#2]}}

\newtheorem{theorem}{Theorem}

\usepackage{algorithmicx}
\usepackage{algorithm}
\usepackage{algpseudocode}
\usepackage{etoolbox}

\makeatletter
\def\blfootnote{\gdef\@thefnmark{}\@footnotetext}
\makeatother

\begin{document}

	\maketitle
	\thispagestyle{empty}
	\pagestyle{empty}
	
	\blfootnote{This work is supported by the grant from the National Robotics Initiative grant titled “NRI: FND:COLLAB: Multi-Vehicle Systems for Collecting Shadow-Free Imagery in Precision Agriculture” (grant no.2019-04791/project accession no. 1020285) from the USDA National Institute of Food and Agriculture.
	}
	
	\maketitle
		\begin{abstract}
			In this work, multiplicative stochasticity is applied to the learning rate of stochastic optimization algorithms, giving rise to stochastic learning-rate schemes. 
			In-expectation theoretical convergence results of Stochastic Gradient Descent equipped with this novel stochastic learning rate scheme under the stochastic setting, as well as convergence results under the online optimization settings are provided. Empirical results consider the case of an adaptively uniformly distributed multiplicative stochasticity and include not only Stochastic Gradient Descent, but also other popular algorithms equipped with a stochastic learning rate. They demonstrate noticeable optimization performance gains, with respect to their deterministic-learning-rate versions. 
		\end{abstract}
		\section{Introduction}
		\label{sec:introduction}
		Machine learning models today range from simple linear regression to deep neural networks. These models try to capture the underlying behavior of systems which produce outputs $ y $ when given inputs $ x $. The performance of these models is usually measured via a loss function $ \ell(h(x;\theta),y) $ where $ h(.;.) $ is the predicted output of the model for an input $ x $ for  given model parameters $ \theta \in \mathbb{R}^d $. Then, if the joint distribution $ P $ of $ x$ and $y $ is known, the estimates $ \theta $ are the minimizers of the expected loss $ \E[x,y\sim P]{\ell(h(x;\theta),y)} $. Nonetheless, complete knowledge of $ P $ is absent in most practical scenarios, and a limited set of training data is available
		instead.
		In that case, an estimate of the expected loss is minimized,
		 known as the empirical loss function, giving rise to the empirical loss minimization (ERM) problem:
		\begin{equation}\label{eq:minimization}
		\mathop { \min }\limits_{\theta  \in {R^d}} \frac{1}{n}\sum\nolimits_{i = 1}^n {\ell(h(x_i;\theta),y_i)} ,
		\end{equation} 
		where the sequence $ \{(x_i,y_i)\}_{i=1}^n $ denotes the available training data.

		Then, if $ \xi $ represents
		a random sample from the
		the available training data,
		let the loss function as a composition of $ \ell(.,.) $ and $ h(.;.) $ be denoted by $L(\theta;\xi)$. Also let $ f(\theta):=\E[\xi]{L(\theta,\xi)} $. This yields a stochastic optimization setting where the computed gradient $ g(\theta,\xi)$ is a function of the true gradient $ \nabla f(\theta) $. If the former is assumed to be an unbiased estimator of the latter as in this work, then: 
		\begin{equation}\label{eq:unbiased_est_grad}
		\mathbb{E}_{\xi_t} [g({\theta_t },{\xi_t })]  = \nabla f(\theta_t),
		\end{equation}
		where $ \xi_t $ and $ \theta_t $ are independent.
		Due to its simplicity, cost-effectiveness and well studied properties, one of the most commonly used algorithms for ERM is stochastic gradient descent (SGD) originally developed in \cite{SGD1951}, which iteratively gives an estimate of $ \theta $ in (\ref{eq:minimization}) via the rule:
		\begin{equation}\label{eq:SGD}
		{\theta _{t + 1}} = {\theta _t} - {a_t}g_t,
		\end{equation}
		for $ t=1,2,\cdots $ and where $ g_t $ is short for $ g({\theta_t },{\xi_t }) $. This rule will be referred to as original SGD or simply SGD throughout the paper. 
	
	Even though commonly the learning rate $ a_t $  is taken to be deterministic, the novelty  is that it introduces the setting where the learning rate in SGD becomes stochastic, by equipping it with multiplicative stochasticity.  That is, the learning rate now becomes $ a_t=\eta_t u_t $  (where e.g., $\eta_t = a/\sqrt{t}$ or $\eta_t=a $, with $a$ a positive constant; with some abuse of terminology, $ \eta_t $ will be referred to as step size) and $ u_t=f(v_b,...,v_t) $ is a function of past random variables $ v_i $, $ i=b,...,t $, referred to as stochastic factors (SF) in the rest of the paper. This work considers two settings, firstly one where the learning rate has memory of all previous SFs, i.e., $ b=1 $, and secondly a setting where the learning rate has no memory of past SFs, i.e., $ b=t $. The latter is analyzed theoretically in the Appendix. Furthermore, a Multiplicative-Stochastic-Learning-Rate is referred to as MSLR. An example of MSLR with memory applied on SGD for  $f(u_1,...,u_t)= \prod\nolimits_{i = 1}^t {{u_i}}$ is shown in Algorithm \ref{alg:MSLR_SGD}.
	
		In general, it has been observed that minimization performance can be noticeably improved under appropriate learning schemes e.g., under deterministic adaptive learning rate schemes with memory (such as memory of past gradients, updates etc). Examples include ADAM \cite{ADAM15} and some variants (e.g., AMSGrad \cite{ADAMGRAD18} or ADAMW \cite{ADAMW2019}) or precursors
		(e.g. \cite{Streeter10}, \cite{Duchi11}). The stochastic learning rate in this work possesses memory by taking into account SFs of previous timesteps to update the parameters of the current timestep. This work demonstrates that performance can be significantly improved if instead of just adaptive, the learning rate becomes both adaptive and stochastic, and possesses memory. In this case, the memory concerns past values of the SF.

		Convergence analysis in expectation but for deterministic learning rates has been studied in \cite{bertsekastsitsiklis2000, nemirovski2009, bachmoulines2011, ghadimi2013,Bottou2018}. This work provides convergence results in-expectation for SGD using stochastic learning rate with memory as well as results under the online learning framework (introduced in \cite{Zinkevich03onlineconvex}). It shows that with stochastic-learning-rates with memory optimization performance is significantly improved, while providing similar technical results to the deterministic-learning-rate case. Concerning memoryless stochastic learning rates, theoretical convergence results in the stochastic approximation and online learning settings are given in the Appendix.
		 Convergence in the almost-surely setting using stochastic learning rate schemes without memory, and in specific MSLR, for algorithms including SGD  has been done in \cite{mamalis2020as_conv_rates}. 
	
\begin{figure}
	\captionsetup[subfigure]{aboveskip=-1pt,belowskip=-5pt}
	\centering
	\begin{subfigure}[b]{\linewidth}
		\centering
		\includegraphics[width=0.5\textwidth]{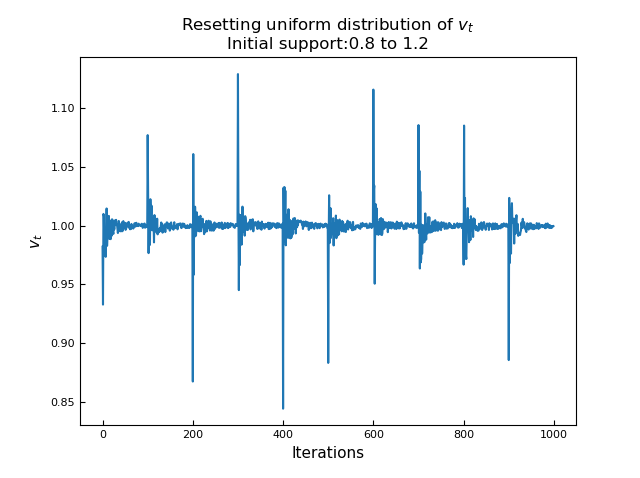}
		\label{fig:lr_0_8_RE}
	\end{subfigure}
	\caption{Plot that illustrates a resetting SF $v_t$, used in the RMSLR-with-memory scheme considered in the experiments, with $v_t \sim \mathcal{U}(\sqrt[\varepsilon_t]{{{c_1}}},\sqrt[\varepsilon_t]{{{c_2}}}) $ where $\varepsilon_t=t \,\, (\bmod \beta)+1$, and $ c_1=0.8, c_2=1.2 $. For the case depicted in the plot, the value of the resetting parameter is $ \beta=100 $, i.e.  the distribution of $ v_t $ resets every $ 100 $ iterations to $ \mathcal{U}(0.8,1.2) $.}
	\label{fig:lr}
\end{figure}

			Moreover, the analysis will allow for resetting SF distributions. Such distributions yield resetting learning rate schemes, which have been shown to improve optimization performance (\cite{odonoghue2012adaptive, loshchilov2017sgdr}). A memoryful Resetting-MSLR scheme will be referred to as RMSLR with memory. A memoryless Resetting-MSLR scheme will be referred to as RMSLR. If, moreover, the  SFs are uniformly distributed, the resetting learning rate scheme will be referred to as RUMSLR with memory or RMSLR for memoryful and memoryless learning rates respectively. If, moreover, its  SFs are uniformly distributed, the resetting learning rate scheme will be referred to as RUMSLR with memory or RUMSLR respectively. An example of a uniformly distributed SF is shown in Fig. \ref{fig:lr}.
		
		The remainder this work is organized as follows. Section \ref{sec:Stochastic Optimization Setting} provides the stochastic optimization setting of this work and provides the in-expectation convergence results. The online learning setting along with the respective convergence results is given in Section \ref{sec:Online Learning Setting}.
		Section \ref{sec:Discussion on the Stochasticity Factor} provides a discussion on the RUMSLR scheme used to obtain the empirical results. 
		
		 In Section \ref{sec:Results}, the experimental results for various optimization algorithms using this memoryful RUMSLR scheme are presented and compared to the deterministic-learning rate case. Section \ref{sec:Conclusion and Future Work} concludes the paper with paths for future work. The Appendix contains results for stochastic learning-rate schemes without memory in the stochastic approximation and online learning rate settings.
			\begin{algorithm}[b]
			\caption{MSLR-with-memory SGD}\label{alg:MSLR_SGD}
			\begin{algorithmic}[]
				\Require
				initial point $ x_1 $, stepsize $ \eta_t $, stochasticity factor $ v_t  $    
				\While{$x_t $ not converged}
				\State \hspace*{-1mm} $t \gets t + 1 $
	
				\State \hspace*{-1mm} $ v_{1:t-1} \gets v_1\cdot...\cdot v_{t-1} $
				\State \hspace*{-1mm} $ a_t \gets \eta_{t} v_{1:t-1} v_k $
				\State \hspace*{-1mm}  $x_{t}\gets x_{t-1} - a_t \nabla f_{v_t}(x_t) $

				\EndWhile\label{alg:beta_RMSLR_SGD_endwhile}
				\State \textbf{return} $x_t$
			\end{algorithmic}
		\end{algorithm}

		\section{Stochastic Optimization Setting}
		\label{sec:Stochastic Optimization Setting}
		This section will provide convergence results for stochastic-learning-rate SGD under the stochastic optimization setting,  briefly discussed in Section \ref{sec:introduction}. These results hold for both MSLR and RMSLR schemes.
		The following assumptions will be used in the theorems. Firstly, function $ f(\theta) $ is $ L $-smooth, i.e.:

		\begin{equation}\label{eq:L_smooth}
		f(\theta ) \le f(\theta ') + \nabla f{(\theta ')^T}(\theta  - \theta ') + \frac{L}{2}{\left\| {\theta  - \theta '} \right\|^2},
		\end{equation}
		for all $ \theta ',\theta  \in \mathbb{R}^d $.
		Furthermore, for all $ \theta \in \mathbb{R}^d$, it satisfies:
		\begin{equation}\label{eq:L_smooth_implication}
		{\left\| {\nabla f(\theta )} \right\|^2} \ge 2c (f(\theta ) - {f_{min}}),
		\end{equation}
		where $ {f_{min}} $ is the global minimum of $ f(\theta) $, and $ 0 < c \le L $. It is noted that this assumption only requires the function to attain a global minimum (more detailed discussion in \cite{Polyak_Convexity2016}).
		Also:
		\begin{equation}\label{eq:variance_related_ineq}
		{\mathbb{E}_{{\xi_t}}}[{\left\| {g({\theta_t},{\xi_t})} \right\|^2}] \le {M} + {M_G}{\left\| {\nabla f({\theta_t})} \right\|^2} \,\, \text{ for } \,\, t \in \mathbb{N},
		\end{equation}
		for some nonnegative scalars $ M $ and $ M_G \ge 1 $.
		Moreover, for the stochastic learning rates it is assumed that the following holds:
\begin{equation} \label{eq:assumption_on_stepsizes}
\begin{aligned}
&\sum\nolimits_{j = 1}^{\infty} {\mathbb{E}_a [a_j]} = \sum\nolimits_{j = 1}^{\infty} \eta_t {\mathbb{E}_a [u_j] } = \infty  
\\ & \sum\nolimits_{j = 1}^{\infty} {\mathbb{V}_a[a_j] } = \sum\nolimits_{j = 1}^{\infty}\eta^2_j {\mathbb{V}_u[u_j]} < \infty
\end{aligned}
\end{equation}
These assumptions are the stochastic-learning-rate equivalents to the typical assumptions made in the deterministic-learning-rate case--commonly referred to as Robbins-Monro conditions (\cite{SGD1951})--namely,  $ \sum\nolimits_j  \eta_j = \infty $ and $ \sum\nolimits_j \eta^2_j < \infty $.
		The theorem that follows provides a bound on the optimality gap of the SGD algorithm when using a stochastic-learning-rate with a constant step size: 
		\begin{theorem}
			\label{th:theorem_finite_dataset_fixed_step_size_R}
			Assume (\ref{eq:L_smooth}-\ref{eq:assumption_on_stepsizes}) are satisfied. Furthermore, assume a positive learning rate of:
			\begin{equation} \label{eq:sqrt_lr_R}
			a_t=  \eta_t u_t , 
			\end{equation}
			where $ \eta_t \in \mathbb{R} $ is either constant or decreasing, and $ u_t := \prod\nolimits_{i = 1}^t {{v_i}} $ where $ v_i $ are i.i.d. 
 with $ \mathbb{E}_v[v_t] \le 1 $, and $ 1  > \lambda \ge \frac{\mathbb{V}_v[v_t]}{\mathbb{E}^2_v[v_t]} $, for some positive $ \lambda \in \mathbb{R} $.
			If, moreover:
			\begin{equation}\label{eq:lr_inequality_fixed_step_size_R_sqrt}
			\eta_t <   \frac{1}{\mathbb{E}_u[u_t] L M_G},
			\end{equation}
			then the iterates of the SGD algorithm in (\ref{eq:SGD}) satisfy:
			\begin{equation} \label{eq:bounds_SGD_constant_step_size_R_sqrt}
			\mathop {\lim }\limits_{t \to \infty } {\mathbb{E} }[f({\theta_{t+1}}) - {f_{min}}]  = 0.
			\end{equation}  
		\end{theorem}
        The guarantee in (\ref{eq:bounds_SGD_constant_step_size_R_sqrt}) is the same as when SGD uses a decreasing deterministic-learning-rate. However, the stepsize's upper bound $\frac{1}{\mathbb{E}_u[u_t] L M_G}$ is larger than the $\frac{1}{L M_G}$ in the deterministic-learning-rate case.   Moreover, the guarantee in (\ref{eq:bounds_SGD_constant_step_size_R_sqrt}) holds for arbitrary decreasing step sizes and does not depend on the initial point unlike decresing stepsizes in previous works (see, e.g., \cite{Bottou2018} for all of the above).         
         In addition, the optimality gap for the stochastic-learning-rate with constant stepsize is zero, whereas for the deterministic-learning-rate case is nonzero, i.e., $\frac{\eta LM}{2c}$. This is also true for a stochastic learning rate scheme without memory in the Appendix
          which yields a nonzero optimality gap for a constant stepsize, whereas Theorem \ref{th:theorem_finite_dataset_fixed_step_size_R}  is able to obtain a zero optimality-gap asymptotically.

		\begin{proof}
			Taking expectations with respect to $ \xi_t $ in (\ref{eq:L_smooth}), and from (\ref{eq:SGD}) substituting $ \theta_{t+1} $ for $ \theta  $ and $ \theta_t $ for $ \theta ' $:
			\begin{equation}
			\begin{aligned}[b] \label{eq:first_expressions_f}
			&{\mathbb{E}_{{\xi_t}}}[f({\theta_{t+1}})]  \le f({\theta_t}) - {a_t}{\nabla ^T}f({\theta_t})\mathbb{E}_{{\xi_t}}[g({\theta_t},{\xi_t})] + a_t^2\frac{L}{2}{\mathbb{E}_{{\xi_t}}}\left[ {{{\left\| {g({\theta_t},{\xi_t})} \right\|}^2}} \right] \\
			& \le   f({\theta_t}) - {a_t}{\left\| {\nabla f({\theta_t})} \right\|^2} + a_t^2\frac{L}{2}{\mathbb{E}_{{\xi_t}}}\left[ {{{\left\| {g({\theta_t},{\xi_t})} \right\|}^2}} \right] \\
			& \le f({\theta_t}) + {a_t}({a_t}\frac{{{M_G}L}}{2} - 1){\left\| {\nabla f({\theta_t})} \right\|^2} + a_t^2\frac{L}{2}M \\
			&  = f({\theta_t}) + \frac{1}{2}{a_t}({a_t}{M_G}L - 2){\left\| {\nabla f({\theta_t})} \right\|^2} + a_t^2\frac{L}{2}M  \\
			&  = f({\theta_t}) + \frac{1}{2}({a^2_t}{M_G}L - 2a_t){\left\| {\nabla f({\theta_t})} \right\|^2} + a_t^2\frac{L}{2}M \\
			& = f({\theta_t}) + \frac{1}{2}(\eta_t^2 {u^2_t}{M_G}L - 2 \eta_t u_t){\left\| {\nabla f({\theta_t})} \right\|^2} + \eta_t^2 u_t^2\frac{L}{2}M .
			\end{aligned}
			\end{equation}
			The second inequality followed from (\ref{eq:unbiased_est_grad}), the third from (\ref{eq:variance_related_ineq}), and the last equality from (\ref{eq:sqrt_lr_R}).			
			Then, taking expectations with respect to $ u_t $ and denoting the result of $ \mathbb{E}_{\xi_t}[\mathbb{E}_{u_t}[.]] $ as $ \mathbb{E}_t[.] $:	
				\begin{equation}
			\begin{aligned}
			&{\mathbb{E}_{{t}}}[f({\theta_{t+1}})]  \le  f({\theta_t}) + \frac{1}{2}({\eta_t^2\mathbb{E}_u[u^2_t]}{M_G}L - 2\eta_t\mathbb{E}_u[u_t]){\left\| {\nabla f({\theta_t})} \right\|^2} + \eta_t^2\mathbb{E}_u[u_t^2]\frac{L}{2}M  	
			\\&  = f({\theta_t}) + \frac{1}{2}\eta_t (\eta_t  \mathbb{E}^2_u[u_t]{M_G}L +  \eta_t \mathbb{V}_u[u_t]{M_G}L - 2 \mathbb{E}_u[u_t]){\left\| {\nabla f({\theta_t})} \right\|^2}  + \eta_t^2\mathbb{E}_u[u_t^2]\frac{L}{2}M ,
			\end{aligned}
			\end{equation}
			where the equality followed from using $ \mathbb{V}_{a}[a_t]=\E[a]{a^2_t}-\E[a]{a_t}^2 $ to replace $ \E[a]{a^2_t} $. Then:	  
			 \begin{equation}
			 \begin{aligned}
			&{\mathbb{E}_{{t}}}[f({\theta_{t+1}})]  = f({\theta_t})  + \frac{1}{2}\eta_t( \mathbb{E}_u[u_t]  +    \frac{\mathbb{V}_u[u_t]}{\mathbb{E}_u[u_t]} - 2\mathbb{E}_u[u_t] ){\left\| {\nabla f({\theta_t})} \right\|^2}  + \eta_t^2\mathbb{E}_u[u_t^2]\frac{L}{2}M  \\
			&  = f({\theta_t}) - \frac{1}{2}\eta_t\mathbb{E}_u[u_t](1-\frac{\mathbb{V}_u[u_t]}{\mathbb{E}^2_u[u_t]} ){\left\| {\nabla f({\theta_t})} \right\|^2}  
			+      \mathbb{E}_a[a^2_t]\frac{L}{2}M 
			\\	&  = f({\theta_t})  - \mathbb{E}_a[a_t]c\varphi_t(f({\theta_t}) - {f_{min}})
			+\mathbb{V}_a[a_t]\frac{L}{2}M   
			+ \mathbb{E}^2_a[a_t]\frac{L}{2}M.
			\end{aligned}	
			\end{equation}
	where $ \varphi_t:=1-\frac{\mathbb{V}_u[u_t]}{\mathbb{E}^2_u[u_t]}  $. This also means that $ \tilde{c} \le \varphi_t <1 $ where $ \tilde{c}=1-\lambda $.	
		Then, adding and subtracting the same terms to create appropriate factor terms:
			\begin{equation}
			\begin{aligned}[b]
			& {\mathbb{E} }[f({\theta_{t+1}})] \le {\mathbb{E} }[f({\theta_t})] - c\varphi_{t}\E[a]{a_t}({\mathbb{E} }[f({\theta_t})] - {f_{min}})	+	 \mathbb{V}_a[a_t]\frac{L}{2}M   
			+ \mathbb{E}^2_a[a_t]\frac{L}{2}M   \\ 
			&+ \E[a]{a_t}\frac{{LM}}{{2c\varphi_{t }}} - \E[a]{a_t}\frac{{LM}}{{2c\varphi_{t  }}}  + \E[a]{a_{t+1}}\frac{{LM}}{{2c\varphi_{t+1}}} - \E[a]{a_{t+1}}\frac{{LM}}{{2c\varphi_{t+1}}}  + {f_{min}} - {f_{min}}.
			\end{aligned}
			\end{equation}
			Then rearranging terms:
			\begin{equation}
			\begin{aligned}[b]
			{\mathbb{E} }[f({\theta_{t+1}}) - {f_{min}}] - \E[a]{a_{t+1}}\frac{{LM}}{{2c\varphi_{t+1}}} &\le (1 - c\varphi_{t }\E[a]{{a_t}})({\mathbb{E} }[f({\theta_t}) - {f_{min}}] - \E[a]{{a_t}}\frac{{LM}}{{2c\varphi_{t}}})
			\\& + (\frac{\E[a]{a_t}}{\varphi_t} - \frac{\E[a]{a_{t+1}}}{\varphi_{t+1}} )\frac{{LM}}{{2c}} + \mathbb{V}_a[{a_t}]\frac{{LM}}{2}. 
			\end{aligned}
			\end{equation}
			By unrolling the sequence it is:
						\begin{equation}
			\begin{aligned}  \label{eq:before_partition_d}
			\mathbb{E}[f({\theta_{t+1}}) - {f_{min}}] - {\mathbb{E}_a}[a_{t+1}]\frac{{LM}}{{2c\varphi_{t+1 }}} &\le \prod\limits_{j=1}^t {(1 - c{\varphi_{j}}{\mathbb{E}_a}[{a_j}])} ({\mathbb{E} }[f({\theta _1}) - {f_{min}}] - \E[a]{a_1}\frac{{LM}}{{2c\varphi_{ 1}}})\\
			&+ \frac{{LM}}{{2c}}\sum\limits_{j = 1}^{t} {( \frac{\E[a]{a_j}}{\varphi_{j}} - \frac{\E[a]{a_{j + 1}}}{\varphi_{j+1}})} \prod\limits_{n = j+1}^{t} {(1 - c{\varphi_{n}}{\mathbb{E}_a}[{a_{n }}])}   \\ 
			&+ \frac{{LM}}{{2}}\sum\limits_{j = 1}^{t} {\mathbb{V}_{a}[{a_j}]} \prod\limits_{n = j+1}^{t} {(1 - c{\varphi_{n}}{\mathbb{E}_a}[{a_{n}}])} .
			\end{aligned}
\end{equation}
Then:
		\begin{equation}
\begin{aligned}  \label{eq:before_partition_d_0}
		\mathbb{E}[f({\theta_{t+1}}) - {f_{min}}] - {\mathbb{E}_a}[a_{t+1}]\frac{{LM}}{{2c\varphi_{t+1 }}} 	&\le  \prod\limits_{j=1}^t {(1 - c{\varphi_{j}}{  \mathbb{E}_a}[{a_j}])} ({\mathbb{E} }[f({\theta _1}) - {f_{min}}] - \E[a]{a_1}\frac{{LM}}{{2c\varphi_{ 1}}})\\
			&+ \frac{{LM}}{{2c}}\sum\limits_{j = 1}^{t} {( \frac{\E[a]{a_j}}{\varphi_{j}} - \frac{\E[a]{a_{j + 1}}}{\varphi_{j+1}})} \prod\limits_{n = j+1}^{t} {(1 - c{\varphi_{n}}{\mathbb{E}_a}[{a_{n }}])}   \\ 
			&+ \frac{{LM}}{{2}}\sum\limits_{j = 1}^{t} { \mathbb{V}_{a}[{a_j}]} \prod\limits_{n = j+1}^{t} {(1 - c{\varphi_{n}}{\mathbb{E}_a}[{a_{n}}])} 
			\\& \le  \prod\limits_{j=1}^t {(1 - c{\varphi_{j}}{  \mathbb{E}_a}[{a_j}])} ({\mathbb{E} }[f({\theta _1}) - {f_{min}}] - \E[a]{a_1}\frac{{LM}}{{2c\varphi_{ 1}}})\\
			&+ \frac{{LM}}{{2c}}\sum\limits_{j = 1}^{t} {( \frac{\E[a]{a_j}}{\varphi_{j}} - \frac{\E[a]{a_{j + 1}}}{\varphi_{j}})} \prod\limits_{n = j+1}^{t} {(1 - c{\varphi_{n}}{\mathbb{E}_a}[{a_{n }}])}   \\  
			&+ \frac{{LM}}{{2}}\sum\limits_{j = 1}^{t} { \mathbb{V}_{a}[{a_j}]} \prod\limits_{n = j+1}^{t} {(1 - c{\varphi_{n}}{\mathbb{E}_a}[{a_{n}}])} ,
			\end{aligned} 
			\end{equation}		
			where the second inequality followed from $- \frac{1}{{{\varphi _j}}} >  - \frac{1}{{{\varphi _{j + 1}}}}$, since ${\varphi _j} > {\varphi _{j + 1}}$, or equivalently, $\frac{{\mathbb{V}[{u_t}]}}{{{E^2}[{u_t}]}} < \frac{{\mathbb{V}[{u_{t + 1}}]}}{{{E^2}[{u_{t + 1}}]}}$. This  is because:
			\begin{equation} \label{eq:variance_of_many_vars}
				\begin{aligned}
				\mathbb{V}[{u_t}] &= \mathbb{V}[\prod\nolimits_{i = 1}^t {{v_i}} ] = \mathbb{E}[\prod\nolimits_{i = 1}^t {v_i^2} ] - \mathbb{E}^2{[\prod\nolimits_{i = 1}^t {{v_i}} ]}= \prod\nolimits_{i = 1}^t {\mathbb{E}[v_i^2]}  - \prod\nolimits_{i = 1}^t {\mathbb{E}^2{{[{v_i}]}}} \\
				&= \prod\nolimits_{i = 1}^t {(\mathbb{V}[{v_i}]}  + \mathbb{E}{[{v_i}]^2}) - \prod\nolimits_{i = 1}^t {\mathbb{E}^2{{[{v_i}]}}}, 
				\end{aligned}
			\end{equation}
since the $v_i$'s are independent. This means that:
\begin{equation} \label{eq:exp_variance_of_many_vars}
\begin{aligned}
&\frac{{\mathbb{V}[{u_t}]}}{{{E^2}[{u_t}]}} = \frac{{\mathbb{V}[\prod\nolimits_{i = 1}^t {{v_i}} ]}}{{{\mathbb{E}^2}[\prod\nolimits_{i = 1}^t {{v_i}} ]}} = \frac{{\prod\nolimits_{i = 1}^t {(\mathbb{V}[{v_i}]}  + \mathbb{E}{{[{v_i}]}^2}) - \prod\nolimits_{i = 1}^t {{\mathbb{E}^2}[{v_i}]} }}{{\prod\nolimits_{i = 1}^t {{\mathbb{E}^2}[{v_i}]} }}= \frac{{\prod\nolimits_{i = 1}^t {(\mathbb{V}[{v_i}]}  + \mathbb{E}{{[{v_i}]}^2})}}{{\prod\nolimits_{i = 1}^t {{\mathbb{E}^2}[{v_i}]} }} - 1.
\end{aligned}
\end{equation}

Using (\ref{eq:exp_variance_of_many_vars}),  the inequality $\frac{{\mathbb{V}[{u_t}]}}{{{\mathbb{E}^2}[{u_t}]}} < \frac{{\mathbb{V}[{u_{t + 1}}]}}{{{\mathbb{E}^2}[{u_{t + 1}}]}}$ can be written as:
\begin{equation}
\begin{aligned}
&\frac{{\prod\nolimits_{i = 1}^t {(\mathbb{V}[{v_i}]}  + \mathbb{E}{{[{v_i}]}^2})}}{{\prod\nolimits_{i = 1}^t {{\mathbb{E}^2}[{v_i}]} }}  < \frac{{\prod\nolimits_{i = 1}^t {(\mathbb{V}[{v_i}] + \mathbb{E}{{[{v_i}]}^2})} (\mathbb{V}[{v_{t + 1}}] + \mathbb{E}{{[{v_{t + 1}}]}^2})}}{{\prod\nolimits_{i = 1}^t {{\mathbb{E}^2}[{v_i}]} {\mathbb{E}^2}[{v_{t + 1}}]}}.
\end{aligned}
\end{equation}
This yields $\frac{{\mathbb{V}[{v_{t + 1}}] + \mathbb{E}^2{{[{v_{t + 1}}]}}}}{{{\mathbb{E}^2}[{v_{t + 1}}]}} > 0$ which holds for all $t$, proving $\frac{{\mathbb{V}[{u_t}]}}{{{\mathbb{E}^2}[{u_t}]}} < \frac{{\mathbb{V}[{u_{t + 1}}]}}{{{\mathbb{E}^2}[{u_{t + 1}}]}}$ and thus $-\frac{1}{{\varphi _{j+1}}} < -\frac{1}{\varphi _{j}}$. Then, from (\ref{eq:lr_inequality_fixed_step_size_R_sqrt}), $ M_G \ge 1 $, and $ c \le L $:
			\begin{equation}
			c \varphi_{t} \E[a]{a_t} \le c   \E[a]{a_t} = c \eta_t \E[u]{u_t}  <   \frac{c }{{ {M_G}L}} \le  \frac{c}{{L}} \le 1,
			\end{equation}
			it is $ 0 \le 1 - c \varphi_{t} \E[a]{a_t}  < 1 $, since $ c \mathbb{E}[a_t] $ is positive, and $ 0<\tilde{c} \le {\varphi_{t}} < 1 $. 

Then, in (\ref{eq:before_partition_d_0}) using $ \tilde{c} \le \varphi_{j} $ and that $ \mathbb{E}_a[a_j] $ is decreasing, it is that:
	\begin{equation}
\begin{aligned}  \label{eq:before_partition_d}
 \mathbb{E}[f({\theta_{t+1}}) - {f_{min}}] - {\mathbb{E}_a}[a_{t+1}]\frac{{LM}}{{2c\varphi_{t+1 }}} &\le  \prod\limits_{j=1}^t {(1 - c{\varphi_{j}}{  \mathbb{E}_a}[{a_j}])} ({\mathbb{E} }[f({\theta _1}) - {f_{min}}] - \E[a]{a_1}\frac{{LM}}{{2c\varphi_{ 1}}})\\
&+ \frac{{LM}}{{2c}}\sum\limits_{j = 1}^{t} {\frac{1}{\tilde{c}}( \E[a]{a_{j}} - \E[a]{a_{j + 1}}  )} \prod\limits_{n = j+1}^{t} {(1 - c{\varphi_{n}}{\mathbb{E}_a}[{a_{n }}])}   \\ 
&+ \frac{{LM}}{{2}}\sum\limits_{j = 1}^{t} { \mathbb{V}_{a}[{a_j}]} \prod\limits_{n = j+1}^{t} {(1 - c{\varphi_{n}}{\mathbb{E}_a}[{a_{n}}])} .
\end{aligned}
\end{equation}

Then, for the second term in (\ref{eq:before_partition_d}) it is that:
			\begin{align}
			& \sum\limits_{j = 1}^{t} \frac{1}{\tilde{c}} {(\E[a]{a_j} - \E[a]{a_{j + 1}})} \prod\limits_{n = j+1}^{t} {(1- c{\mathbb{E}_a}[{a_{n }}])}   \le \sum\limits_{j = 1}^{t} \frac{1}{\tilde{c}} {(\E[a]{a_j} - \E[a]{a_{j + 1}})}   \nonumber =\frac{1}{\tilde{c}} (\E[a]{a_1} - \E[a]{a_{t + 1}}),
			\end{align} 
			which converges to $\frac{1}{\tilde{c}} \E[a]{a_1} $ as $ t \to \infty $ since $ \lim\limits_{t \to \infty } \E[a]{a_{t + 1}} = 0 $. 
			For the third term in (\ref{eq:before_partition_d}) it is that:
			\begin{equation}\label{eq:before_partition}
\begin{aligned}[b]
&\sum\limits_{j = 1}^{t}{\mathbb{V}_{a}[{a_j}]} \prod\limits_{n = j+1}^{t} {(1 - \varphi_n c{\mathbb{E}_a}[{a_{n}}])}  \le \sum\limits_{j = 1}^{t}{\mathbb{V}_{a}[{a_j}]} < \infty,
\end{aligned}
\end{equation}
from the assumptions of the theorem.
			Then let the algorithm run for long enough, so that the last two terms in (\ref{eq:before_partition_d}) are close enough to their limits so that for their tails after some timestep $ K $ the following holds:
			\begin{equation}
			\begin{aligned}[b] \label{eq:aprox_eq_to_zero_d}
			\sum\limits_{j = K+1}^{t} \frac{1}{\tilde{c}}{(  \E[a]{a_j} - \E[a]{a_{j + 1}}  ) } \prod\limits_{n = j+1}^{t} {(1 - c\varphi_n{\mathbb{E}_a}[{a_{n }}])} & \approx 0 \\    
			\sum\limits_{j = K+1}^{t}{\mathbb{V}_{a}[{a_j}]} \prod\limits_{n = j+1}^{t} {(1 - c\varphi_n  {\mathbb{E}_a}[{a_{n}}])} & \approx 0.
			\end{aligned}
			\end{equation}
			Then, partitioning the last two terms in (\ref{eq:before_partition_d}) at $ K $, and using (\ref{eq:aprox_eq_to_zero_d}), it is that:
			\begin{equation}
			\begin{aligned}[b] \label{eq:partition_d}
			\mathbb{E}[f({\theta_{t+1}}) - {f_{min}}] - {\mathbb{E}_a}[a_{t+1}]\frac{{LM}}{{2c\varphi_{t+1}}} &\le \prod\limits_{j=1}^t {(1 - c\varphi_{j}{\mathbb{E}_a}[{a_j}])} ({\mathbb{E} }[f({\theta _1}) - {f_{min}}] - \E[a]{a_1}\frac{{LM}}{{2c\varphi_{ 1}}}) 
			\\	& + \frac{{LM}}{{2c}}\sum\limits_{j = 1}^{K} \frac{1}{\tilde{c}} {(  \E[a]{a_j} - \E[a]{a_{j + 1}} )  } \prod\limits_{n = j+1}^{t} {(1 - c\varphi_{n}{\mathbb{E}_a}[{a_{n }}])}   
			\\	& + \frac{{LM}}{{2}}\sum\limits_{j = 1}^{K} {\mathbb{V}_{a}[{a_j}]} \prod\limits_{n = j+1}^{t} {(1 - c\varphi_{n}{\mathbb{E}_a}[{a_{n}}])}.
			\end{aligned}
			\end{equation}
			Therefore, from (\ref{eq:partition_d}), and from $ \E[a]{a_{t+1}}=\eta_{t+1}\E{u_{t+1}} $ it follows that:
			\begin{equation} \label{eq:decreasing_step_size}
			\mathop {\lim }\limits_{t \to \infty } \left(	{\mathbb{E} }[f({\theta_{t+1}}) - {f_{min}}] - 
			\eta_{t+1}\E{u_{t+1}}\frac{{LM}}{{2c}} \right) = 0,
			\end{equation}
			which means:
			\begin{equation}
			\mathop {\lim }\limits_{t \to \infty } {\mathbb{E} }[f({\theta_{t+1}}) - {f_{min}}]  =0,
			\end{equation}
			since $  \mathop {\lim }\limits_{t \to \infty }\eta_{t+1}\E{u_{t+1}}\frac{{LM}}{{2c}} =0$ from $    \mathop {\lim }\limits_{t \to \infty } \E{u_{t+1}}  =0 $ and $ \eta_{t+1} $ being either constant or decreasing.
			%
		\end{proof}

Concluding this section, from Theorem 
\ref{th:theorem_finite_dataset_fixed_step_size_R} it is thus observed that for the stochastic optimization setting the optimality gap bounds for SGD using a stochastic-learning-rate are improved from the ones when using a deterministic-learning-rate. This is because, the optimality gap is zero regradless of whether a constant or a decreasing learning rate is used, whereas for a deterministic learning rate, the optimality gap when using a determinisitc learning rate is nonzero (e.g. see \cite{Bottou2018} for more details).
		\section{Online Learning Setting}
		\label{sec:Online Learning Setting}
		Now, an alternative analysis for stochastic-learning-rate SGD which holds for both MSLR and RMSLR is given by considering the following optimization problem, referred to as online optimization:
		\begin{equation}\label{eq:online_minimization}
		{R_T} = \sum\limits_{t = 1}^T {{f_t}({\theta _t})}  - \mathop {\min }\limits_{\theta  \in X} \sum\limits_{t = 1}^T {{f_t}({\theta })} ,
		\end{equation}
		where $ {\theta ^*} = \mathop {\arg \min }_{\theta  \in X} \sum\nolimits_{t = 1}^T {{f_t}(\theta )} $, functions $ f_1$ through $f_T $ are convex, $ g_t:=\nabla f_t(\theta_t) $, and $ X \subseteq {R^d} $ is a convex and closed set. The function $ R_T $ is called the regret up to timestep $ T $. This optimization framework was introduced in \cite{Zinkevich03onlineconvex}. Then, let the online algorithm:
		\begin{equation}\label{eq:online_SGD}
		\begin{aligned}[b]
		{ \nu _{t}} = {\theta _t} - {a_t}g_t \\
		{\theta _{t + 1}} = \Pi _X ({\nu _t}),
		\end{aligned}
		\end{equation}
		where $ {\Pi _X}(y) = \mathop {\arg \min } _{\theta  \in X} \left\| {x - y} \right\| $ is the projection of $ y $ on $ X $.  Studying an online algorithm with vanishing average regrets amounts to studying a stochastic algorithm for ERM as discussed in \cite{Bianchi04}.
		
		The following assumptions are used in the online minimization setting. The first is that functions $ f_t(\theta) $ are convex, i.e.:
		\begin{equation}\label{eq:f_convex}
		f_t(\theta ) \ge f_t(\theta ') + \nabla f_t{(\theta ')^T}(\theta  - \theta ') \,\, \text{ for all } \,\,  \theta ',\theta  \in \mathbb{R}^d.
		\end{equation}
		Furthermore, if $ X $ is a set with bounded diameter $ D $ then:
		\begin{equation}
		\left\| {{\theta'} - \theta} \right\| \le D ,
		\end{equation}
		for all $ {{\theta'}, \theta} \in X$.
		Moreover, if $X \subseteq {R^d}$ is a convex and closed set, then for any $ \theta \in X $ and $ y \in \mathbb{R}^d $, it is that:
		\begin{equation}\label{eq:projection_property}
		\left\| {y - \theta} \right\| \ge \left\| {{\theta'} - \theta} \right\|,
		\end{equation}
		for	$ {\theta'} = {\Pi _X}(y) = \mathop {\arg \min }_{\theta  \in X} \left\| {y - \theta} \right\| $. 
		It is noted, that the regret bounds for the deterministic-learning-rate version of SGD are obtained by following the proofs of the next Theorem  \ref{th:theorem_online_R} but using $ a_t= \frac{a}{\sqrt{t}} $ as the learning rate.

		Then, it will be shown that stochastic-learning-rate SGD has vanishing average regret for diminishing step size $ \eta_t = \frac{a}{\sqrt{t}} $.  The theorem for the online learning setting follows:
				\begin{theorem}  
			\label{th:theorem_online_R}
			Assume that $ X $ is convex, closed, and has bounded diameter $ D $, and that the gradient of $ f_t $ is bounded, i.e. $ \left\| {\nabla {f_t}(\theta )} \right\|_\infty \le G $. Furthermore, assume a positive learning rate of
			$ 	a_t=\frac{a}{\sqrt{t}} u_t $
			with $ a = \frac{D}{{2G }} $.	Then 
			for all $ T\ge 1 $ it is that:
			\begin{equation} \label{eq:online_mslrsgd_beta}
			{R_{a,T}} \le 2 DG\sqrt T ,
			\end{equation}
		\end{theorem}
		where the subscript $ a $ denotes the expectation of $ R_T $ with respect to the stochastic learning rate sequence $\{a_t\}_{t=1}^{T} $. Therefore, $ R_{a,T}/T \to 0 $ for $ T \to \infty $.
	
		\begin{proof}
			Firstly, it is noted that:
			\begin{equation}
			\begin{aligned}[b]
			{\left\| {{\theta _{t + 1}} - {\theta ^*}} \right\|^2} &= {\left\| {{\Pi _X}({\theta _t} - {a_t}{g_t}) - {\theta ^*}} \right\|^2}\\
			&\le {\left\| {{\theta _t} - {a_t}{g_t} - {\theta ^*}} \right\|^2} \\
			&= {\left\| {{\theta _t} - {\theta ^*}} \right\|^2} - 2{a_t}g_t^T({\theta _t} - {\theta ^*}) + a_t^2{\left\| {{g_t}} \right\|^2},
			\end{aligned}
			\end{equation}
			where the first inequality followed from (\ref{eq:projection_property}) replacing $\theta $ by $ \theta^* $ and $ y $ by $ {\theta _t} - {a_t}{g_t} $.
			Taking expectations with respect to the learning rates  $ a_1\cdots a_t $ and denoting the results as $  \E[{a}]{.} $,  it is that:
			\begin{equation}
			\begin{aligned}[b]
			&\mathbb{E}_a[{\left\| {{\theta _{t + 1}} - {\theta ^*}} \right\|^2}] 
		 \le   \mathbb{E}_a[{\left\| {{\theta _t} - {\theta ^*}} \right\|^2}]   
		 \\&- 2\E[{a}]{a_t} \mathbb{E}_a [g_t^T({\theta _t} - {\theta ^*})] + \E[{a}]{a_t^2}\mathbb{E}_a[{\left\| {{g_t}} \right\|^2}],
			\end{aligned}
			\end{equation}
			 Rearranging with respect to $  g_t^T({\theta _t} - {\theta ^*}) $ and using (\ref{eq:f_convex}), it is:
			\begin{equation}
			\begin{aligned}[b]
			&\mathbb{E}_a[{f_t}({\theta _t}) - {f_t}({\theta ^*})] \le \mathbb{E}_a[g_t^T({\theta _t} - {\theta ^*})]\\
			&\le \frac{{{ \mathbb{E}_a[{\left\| {{\theta _t} - {\theta ^*}} \right\|}^2]} - \mathbb{E}_{a}[{{\left\| {{\theta _{t + 1}} - {\theta ^*}} \right\|}^2}] }}{{2{\E[{a}]{a_t}}}} + \frac{{{\E[a]{a_t^2}}}}{2\E[a]{a_t}}\mathbb{E}_a[{\left\| {{g_t}} \right\|^2}]
			\\& \le \frac{{\mathbb{E}_a[{{\left\| {{\theta _t} - {\theta ^*}} \right\|}^2}] - \mathbb{E}_{a}[{{\left\| {{\theta _{t + 1}} - {\theta ^*}} \right\|}^2}] }}{{2{\E[a]{a_t}}}} + \frac{{{\E[a]{a_t} }}}{2} \mathbb{E}_a[{\left\| {{g_t}} \right\|^2}]
			\\& +\frac{{{\mathbb{V}_{a}[{a_t}]}}}{2{\E[a]{a_t}}}\mathbb{E}_a[{\left\| {{g_t}} \right\|^2}].
			\end{aligned}
			\end{equation}
			Then, summing over all steps:
\begin{equation} \label{eq:summing_all}
\begin{aligned}[b]
\sum\limits_{t = 1}^T \mathbb{E}_a[{[{f_t}(\theta ) - {f_t}({\theta ^*})]}]  & \le \sum\limits_{t = 1}^T {\frac{1}{{2{\mathbb{E}_{a}[a_t]}}}\left( {\mathbb{E}_a[{{\left\| {{\theta _t} - {\theta ^*}} \right\|}^2}] - \mathbb{E}_{a} [{{\left\| {{\theta _{t + 1}} - {\theta ^*}} \right\|}^2}] } \right)}\\ & + \frac{1}{2}\sum\limits_{t = 1}^T {\E[a]{a_t} } \mathbb{E}_a[{\left\| {{g_t}} \right\|^2}] + \frac{1}{2}\sum\limits_{t = 1}^T \frac{{\mathbb{V}_{a}[{a_t}]}}{\E[a]{a_t}}{ \mathbb{E}_a[{{\left\| {{g_t}} \right\|}^2}]} 
\\ &= \frac{1}{{2{\mathbb{E}_a[a_1]}}}\mathbb{E}_a[{\left\| {{\theta _1} - {\theta ^*}} \right\|^2}] 
\\& + \sum\limits_{t = 2}^T { \mathbb{E}_{a}{{\left\| {{\theta _t} - {\theta ^*}} \right\|}^2}\left( {\frac{1}{{2{\E[a]{a_t}}}} - \frac{1}{{2\E[a]{a_{t-1}} }}} \right)}\\& + \frac{1}{2}\sum\limits_{t = 1}^T \E[a]{a_t}{\mathbb{E}_a[\left\| {{g_t}} \right\|^2]} + \frac{1}{2}\sum\limits_{t = 1}^T \frac{ {\mathbb{V}_{a}[{a_t}]}}{ \E[a]{a_t}}{ \mathbb{E}_a[{{\left\| {{g_t}} \right\|}^2}] }.
\end{aligned}
\end{equation} 
Using the diameter of the convex set $ D $:
\begin{equation} \label{eq:summing_all}
\begin{aligned}[b]
 \sum\limits_{t = 1}^T \mathbb{E}_a[{[{f_t}(\theta ) - {f_t}({\theta ^*})]}] & \le \frac{1}{{2{\mathbb{E}_a[a_1]}}}D^2 
 + \sum\limits_{t = 2}^T { D^2\left( {\frac{1}{{2{\E[a]{a_t}}}} - \frac{1}{{2\E[a]{a_{t-1}} }}} \right)}\\& + \frac{1}{2}\sum\limits_{t = 1}^T \E[a]{a_t}{\mathbb{E}_a[\left\| {{g_t}} \right\|^2]} + \frac{1}{2}\sum\limits_{t = 1}^T \frac{ {\mathbb{V}_{a}[{a_t}]}}{ \E[a]{a_t}}{ \mathbb{E}_a[{{\left\| {{g_t}} \right\|}^2}] }.
\end{aligned}
\end{equation}
Expanding the telescoping sum:
\begin{equation} \label{eq:summing_all}
\begin{aligned}[b]
\sum\limits_{t = 1}^T \mathbb{E}_a[{[{f_t}(\theta ) - {f_t}({\theta ^*})]}] & \le \frac{1}{{2{\mathbb{E}_a[a_1]}}} D^2 
+ { D^2 \left( {\frac{1}{{2{\E[a]{a_T}}}} - \frac{1}{{2\E[a]{a_{1}} }}} \right)} \\& + \frac{1}{2}\sum\limits_{t = 1}^T \E[a]{a_t}{\mathbb{E}_a[\left\| {{g_t}} \right\|^2]} + \frac{1}{2}\sum\limits_{t = 1}^T \frac{ {\mathbb{V}_{a}[{a_t}]}}{ \E[a]{a_t}}{ \mathbb{E}_a[{{\left\| {{g_t}} \right\|}^2}] }.
\end{aligned}
\end{equation}

This means that:
\begin{equation} \label{eq:summing_all}
\begin{aligned}[b]
\sum\limits_{t = 1}^T \mathbb{E}_a[{[{f_t}(\theta ) - {f_t}({\theta ^*})]}] \le \frac{1}{{2{\mathbb{E}_a[a_T]}}} D^2 
 + \frac{1}{2}\sum\limits_{t = 1}^T \E[a]{a_t}{\mathbb{E}_a[\left\| {{g_t}} \right\|^2]} + \frac{1}{2}\sum\limits_{t = 1}^T \frac{ {\mathbb{V}_{a}[{a_t}]}}{ \E[a]{a_t}}{ \mathbb{E}_a[{{\left\| {{g_t}} \right\|}^2}] },
\end{aligned}
\end{equation}
where the third line followed by using the diameter $ D $ of the convex set $ X $, and the fourth line by using the telescoping sum of the middle term.
			Using $ a_t= \frac{a}{\sqrt{t}} u_t $ it is:
						\begin{equation}
			\begin{aligned}[b]
			\sum\limits_{t = 1}^T { \mathbb{E}_a[[{f_t}(\theta ) - {f_t}({\theta ^*})]] }  \le \frac{1}{{2{\mathbb{E}_a[a_T]}}}D^2
			 + \frac{1}{2}\sum\limits_{t = 1}^T \frac{a}{\sqrt{t}}\mathbb{E}_u[u_t]{\mathbb{E}_a[\left\| {{g_t}} \right\|^2]} + \frac{1}{2}\sum\limits_{t = 1}^T \frac{a{\mathbb{V}_{u}[{u_t}]}}{\sqrt{t}\E[u]{u_t}^2}{ \mathbb{E}_a[{{\left\| {{g_t}} \right\|}^2}] }.			\end{aligned}
			\end{equation}
This means that, using $ \mathbb{E}_u[u_t] \le 1 $ and $\frac{\mathbb{V}_{u}[{u_t}]}{ \mathbb{E}_u^2[u_t]}  \le 1$:
			\begin{equation}
\begin{aligned}[b]
\sum\limits_{t = 1}^T { \mathbb{E}_a[[{f_t}(\theta ) - {f_t}({\theta ^*})]] }   
& \le \frac{1}{{2{\mathbb{E}_a[a_T]}}} D^2
+ \frac{1}{2}\sum\limits_{t = 1}^T \frac{a}{\sqrt{t}}{\mathbb{E}_a[\left\| {{g_t}} \right\|^2]} + \frac{1}{2}\sum\limits_{t = 1}^T \frac{a}{\sqrt{t}}{ \mathbb{E}_a[{{\left\| {{g_t}} \right\|}^2}] }
\\ &\le \frac{\sqrt{T}}{{2a}} D^2  + \sum\limits_{t = 1}^T \frac{a}{\sqrt{t}}{ \mathbb{E}_a[{{\left\| {{g_t}} \right\|}^2}] }
\le \frac{{{D^2} \sqrt{T} }}{{2a}} +  aG^2\sum\limits_{t = 1}^T {\frac{1}{{\sqrt t }}}
\\& \le \frac{{{D^2} \sqrt{T} }}{{2a }} +  aG^2 2\sqrt{T}.
\end{aligned}
\end{equation}
			Then, minimizing the right hand side of the third inequality yields $ a ={D}/({{2G  }})$ for which:
			\begin{equation}
			{R_{a,T}} \le   2DG\sqrt T ,
			\end{equation}
			thus concluding the proof.
		\end{proof}
		In conclusion, from Theorem 
		\ref{th:theorem_online_R} it is thus observed that for the online optimization setting the expected--with respect to the learning rate--regret bound yield the same expected regret rate $ O(\frac{1}{\sqrt{T}}) $  as the deterministic-learning-rate SGD algorithm (e.g., in \cite{Duchi11}). 
		\section{Discussion on the Stochasticity Factor}
		\label{sec:Discussion on the Stochasticity Factor}
		In this work, the stochastic learning rate scheme considered for the empirical results is a RMSLR-with-memory scheme and specifically $ u_t=\prod_{i=1}^{t}v_t $ with $v_t \sim \mathcal{U}(\sqrt[\varepsilon_t]{{{c_1}}},\sqrt[\varepsilon_t]{{{c_2}}}) $ where $\varepsilon_t=t \,\, (\bmod \beta)+1$, and $ \beta \ge 1 $, shown in Fig. \ref{fig:lr}. In this case the distribution of $ v_t $ resets every $ \beta $ number of iterations to $ \mathcal{U}(c_1,c_2) $, named $ \beta-$RUMSLR-with-memory because the SF is uniformly distributed.
		Moreover, it has been observed that choosing low values for $ c_1 $ or large values for $ c_2 $ slows down the algorithm significantly or destabilizes it. On the other side, they should be kept adequately apart so that the bursts are significant enough.

		The distribution of $ v_t $ is chosen to be reset after it is adequately close to unity and enough, $ \beta $ timesteps have passed, e.g., as shown in Fig. \ref{fig:lr}, where $ \beta =100$, indicated by the bursts that occur every 100 iterations.  Furthermore, it can be shown that the condition of Theorem \ref{th:theorem_finite_dataset_fixed_step_size_R} that $  \frac{{{\mathbb{V}_u}[{u_t}]}}{{\mathbb{E}_u^2[{u_t}]}} \le \lambda $ holds with $ \lambda=\frac{1}{3} $ for all $ c_1,c_2 $. This is because $ {(\sqrt[{{\varepsilon _t}}]{{{c_2}}} - \sqrt[{{\varepsilon _t}}]{{{c_1}}})^2} < {(\sqrt[{{\varepsilon _t}}]{{{c_2}}} + \sqrt[{{\varepsilon _t}}]{{{c_1}}})^2} $ which means that $ \frac{{{V_v}[{v_t}]}}{{{\mathbb{E}_v}[{v_t}]}}= \frac{{4{{(\sqrt[{{\varepsilon _t}}]{{{c_2}}} - \sqrt[{{\varepsilon _t}}]{{{c_1}}})}^2}}}{{12{{(\sqrt[{{\varepsilon _t}}]{{{c_2}}} + \sqrt[{{\varepsilon _t}}]{{{c_1}}})}^2}}}  \le \frac{1}{3}$, since $ v_t $ is uniformly distributed. This means that:
\begin{equation}
\begin{aligned}
&\frac{{{\mathbb{V}_u}[{u_t}]}}{{\mathbb{E}_u^2[{u_t}]}}=\frac{{\prod\nolimits_{i = 1}^t {({\mathbb{V}_u}[{u_t}] + {\mathbb{E}_u}[{u_t}])} }}{{\prod\nolimits_{i = 1}^t {{\mathbb{E}_u}[{u_t}]} }} - 1 = \prod\nolimits_{i = 1}^t {\left( {\frac{{{\mathbb{V}_v}[{v_t}]}}{{{\mathbb{E}_v}[{v_t}]}} + 1} \right)}  - 1 
\le {\left( {\frac{1}{3} + 1} \right)^t} - 1,
\end{aligned}
\end{equation} 
which yields $ \frac{{{\mathbb{V}_u}[{u_t}]}}{{\mathbb{E}_u^2[{u_t}] }} \le \frac{1}{3} $, for $ t=1 $, which also holds for all $ t $ since $ {\left( {\frac{1}{3} + 1} \right)}^t - 1 $ is increasing. Additionally, it can be checked numerically that the condition $ \mathbb{E}_u^2[{v_t}] \le 1 $
as well as that the stochastic Robins-Monro conditions in Assumption \ref{eq:assumption_on_stepsizes} are satisfied for the RUMSLR scheme used to derive the empirical results in Section (\ref{sec:Results}). It is noted that, since the aforementioned conditions are sufficient and not necessary,  stochastic learning rates that do not satisfy these conditions but still improve performance might exist. 
		
		Experimentally, it will be demonstrated that  $ \beta-$RUMSLR-with-memory SGD yields noticeable improvements in minimization performance for datasets such as CIFAR-10, and CIFAR-100 (\cite{CIFARdata09}) when compared to SGD without stochastic learning rate. This scheme was equipped to the learning rates of other algorithms as well such as ADAM \cite{ADAM15}, AMSGrad \cite{ADAMGRAD18}, ADAMW \cite{ADAMW2019}, and two SGD with momentum algorithms 
		originally presented in \cite{SGDM1964} and \cite{ASGDM1983} (subsequently referred to as SGD with momentum first and second versions respectively),
		all of which, even though not theoretically studied in this work, yielded similar improved performance when using $ \beta-$RUMSLR with memory.

		\section{Results}
		\label{sec:Results}
		Experimental results using the $ \beta-$RUMSLR-with-memory  scheme described in \ref{sec:Discussion on the Stochasticity Factor} applied to SGD, SGD with momentum versions first and second,  ADAM, AMSGrad, and ADAMW are presented in this section. For $ \beta-$RUMSLR-with-memory with memory the plots show the average running loss every 20 mini-batches for CIFAR-10 and CIFAR-100, and every 200 mini-batches for MNIST.  
		The batch size was 128 in all experiments. The rest of the hyperparameters, i.e., learning rates, momentums (when applicable) and SF constants, are chosen via grid search.
		The datasets were standardized before use.
		For CIFAR-10 and CIFAR-100, the network architecture was Pytorch's ResNet18.
		For MNIST, logistic regression was used. For the former two datasets a constant stepsize was used, and a diminishing stepsize for the latter, for all algorithms.

\begin{figure}
			\centering
			\begin{subfigure}[b]{0.44\textwidth}
				\centering
				\includegraphics[width=\textwidth]{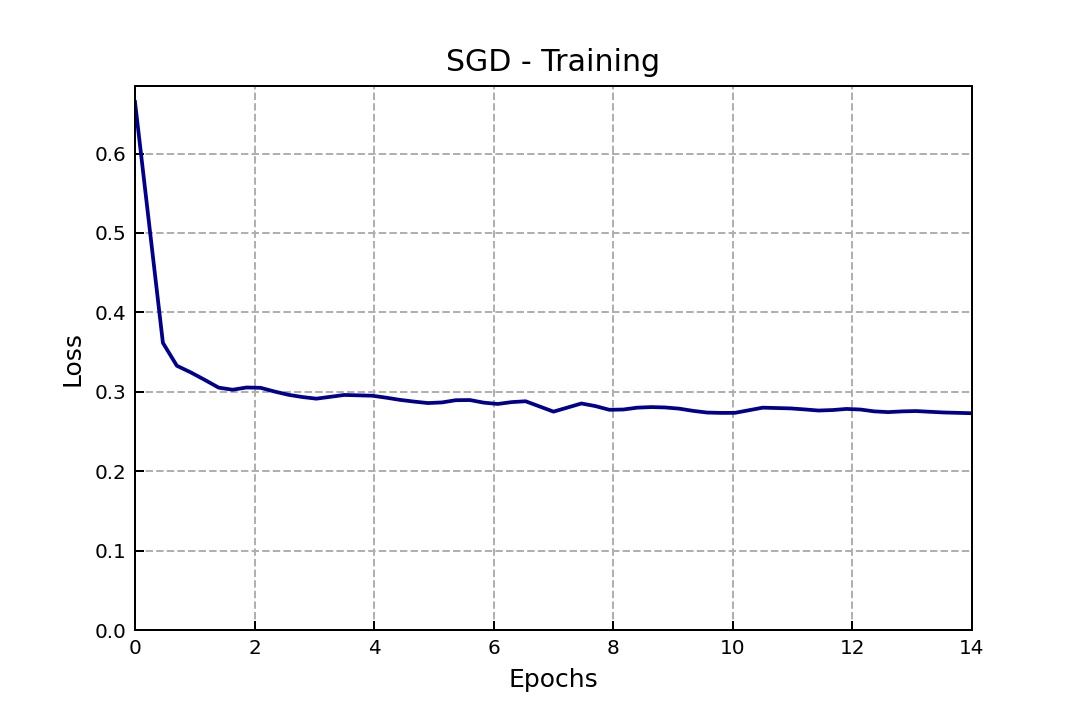}
				\label{fig:}
			\end{subfigure} %
			\begin{subfigure}[b]{0.44\textwidth}
				\centering
				\includegraphics[width=\textwidth]{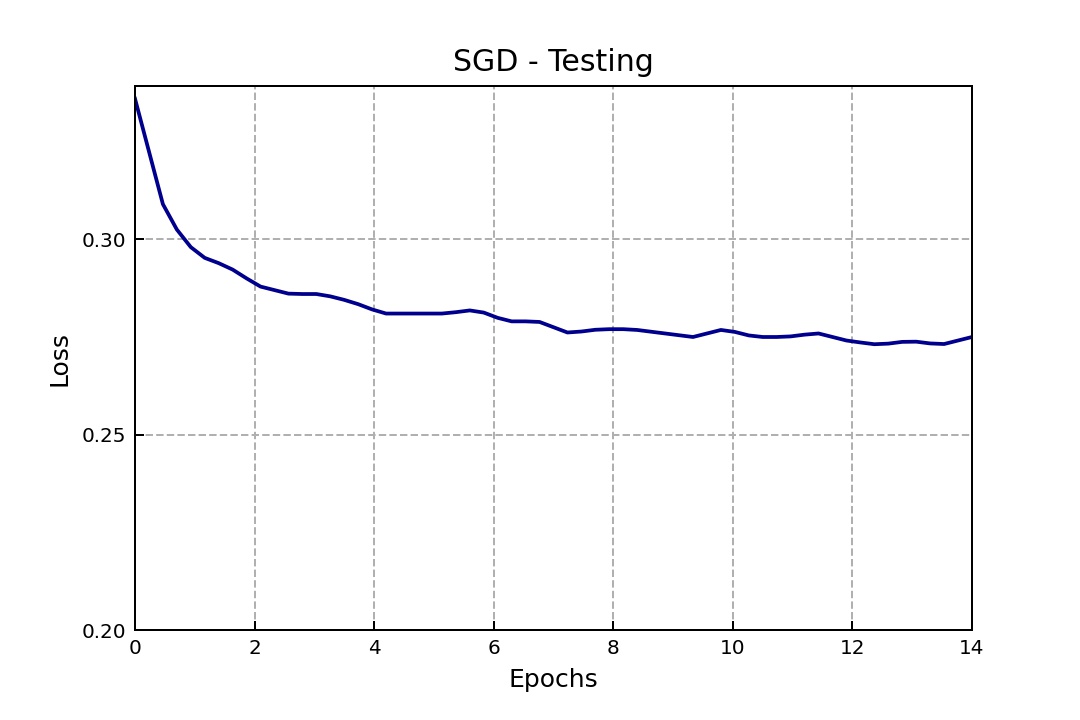}
				\label{fig:}
			\end{subfigure} %
			\caption{Plots that illustrate convergence for the $\beta$-RUMSLR-with-memory SGD for $ \beta=100 $, $ c_1=0.7 $ and $ c_2=1.3 $ in a convex setting.}
			\label{fig:MNIST}
		\end{figure}

\begin{figure}			
			\centering
			\begin{subfigure}[b]{0.32\textwidth}
				\centering
				\includegraphics[width=\textwidth]{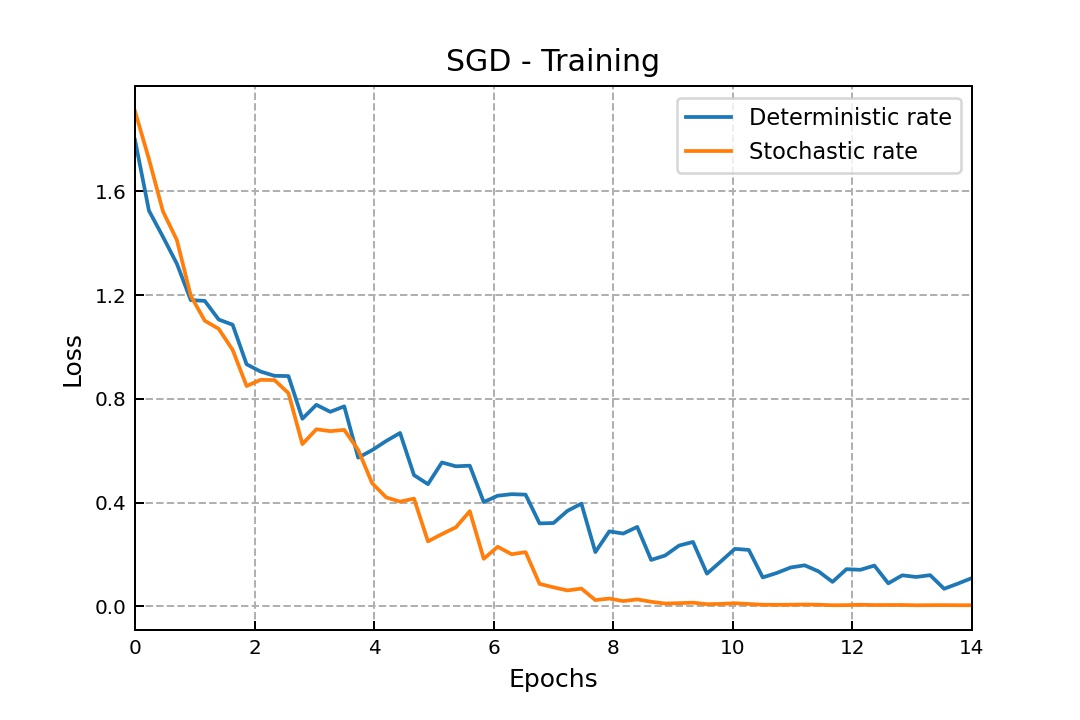}
				\label{fig:}
			\end{subfigure}
			\begin{subfigure}[b]{0.32\textwidth}
				\centering
				\includegraphics[width=\textwidth]{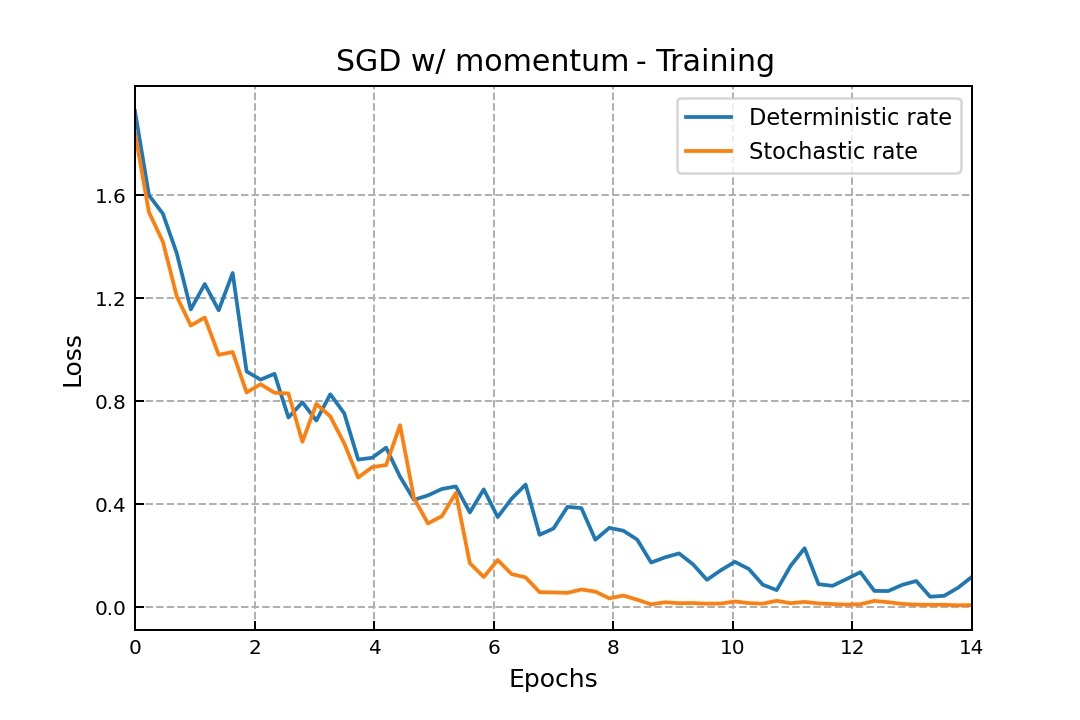}
				\label{fig:}
			\end{subfigure}
			\begin{subfigure}[b]{0.32\textwidth}
				\centering
				\includegraphics[width=\textwidth]{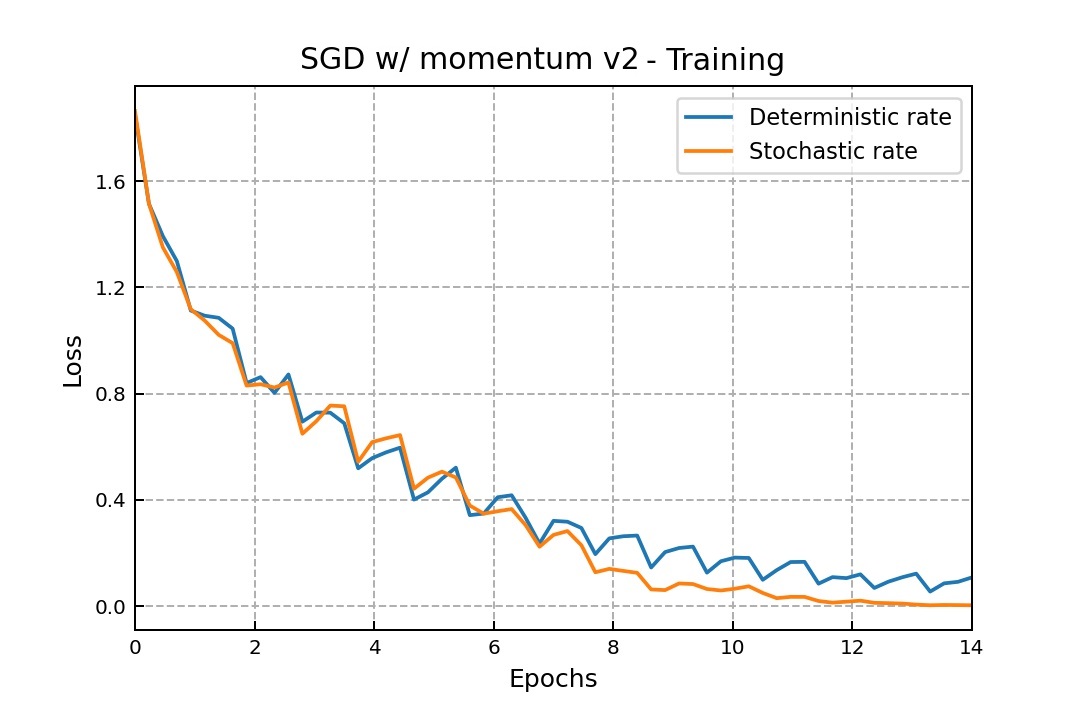}
				\label{fig:}
			\end{subfigure}
			\begin{subfigure}[b]{0.32\textwidth}
				\centering
				\includegraphics[width=\textwidth]{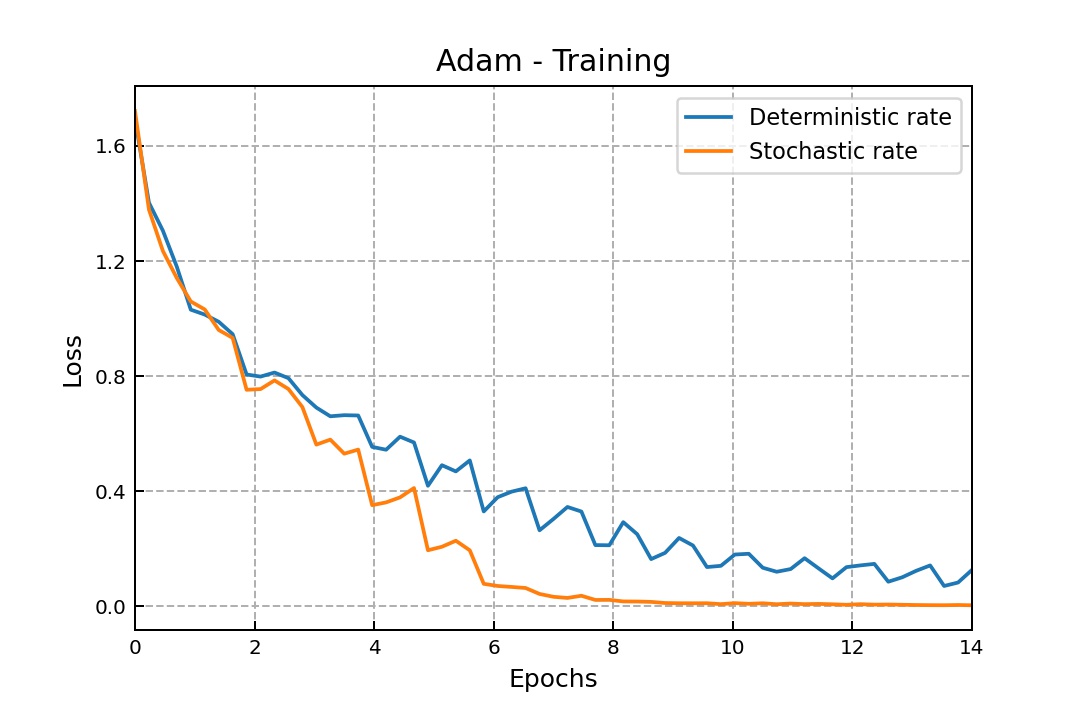}
				\label{fig:}
			\end{subfigure}
			\begin{subfigure}[b]{0.32\textwidth}
				\centering
				\includegraphics[width=\textwidth]{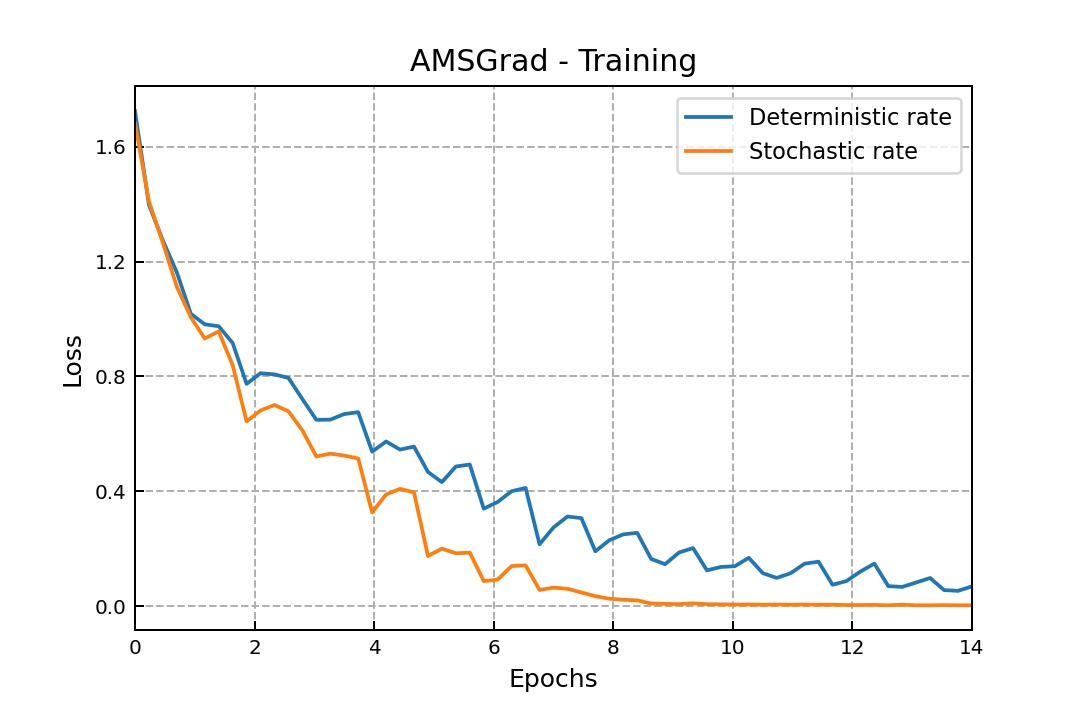}
				\label{fig:}
			\end{subfigure}
			\begin{subfigure}[b]{0.32\textwidth}
				\centering
				\includegraphics[width=\textwidth]{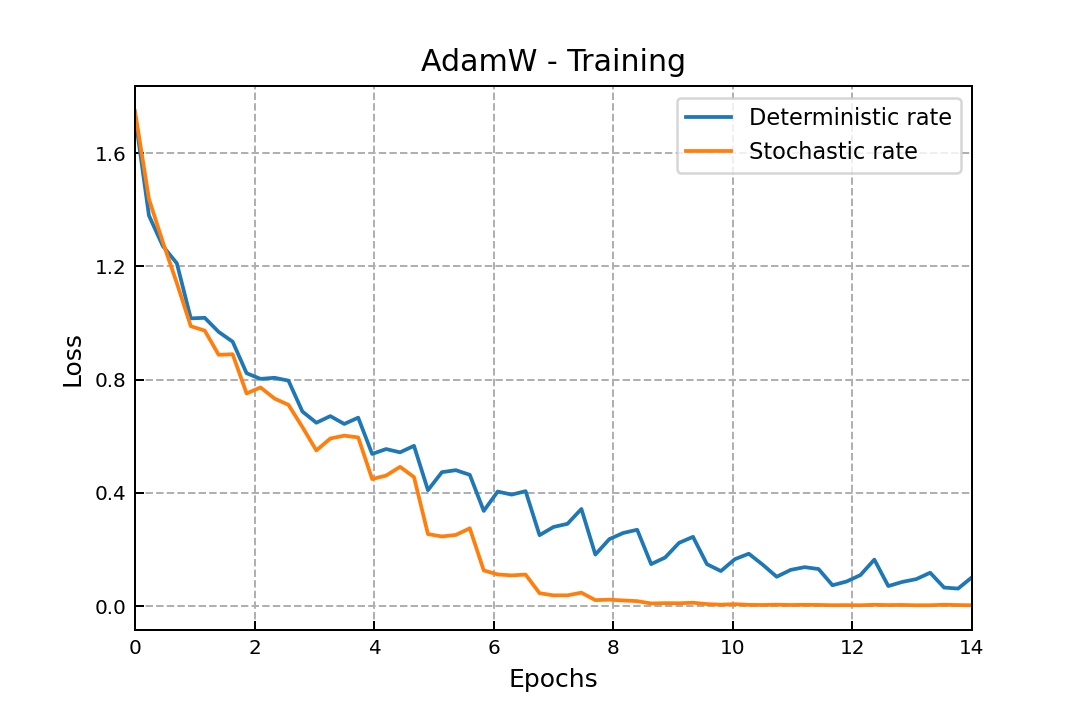}
				\label{fig:}
			\end{subfigure}
			\vspace{-5mm}
			\caption{Plots that illustrate the  $\beta$-RUMSLR-with-memory scheme (in orange) and the deterministic learning rate (in blue) schemes for $ \beta=100 $, $ c_1=0.8 $ and $ c_2=1.2 $. It is observed that for the CIFAR-10 dataset algorithms with  $\beta$-RUMSLR yield noticeable gains in minimization performance.}
			\label{fig:CIFAR10}
		\end{figure}

\begin{figure}
			\centering
			\begin{subfigure}[b]{0.32\textwidth}
				\centering
				\includegraphics[width=\textwidth]{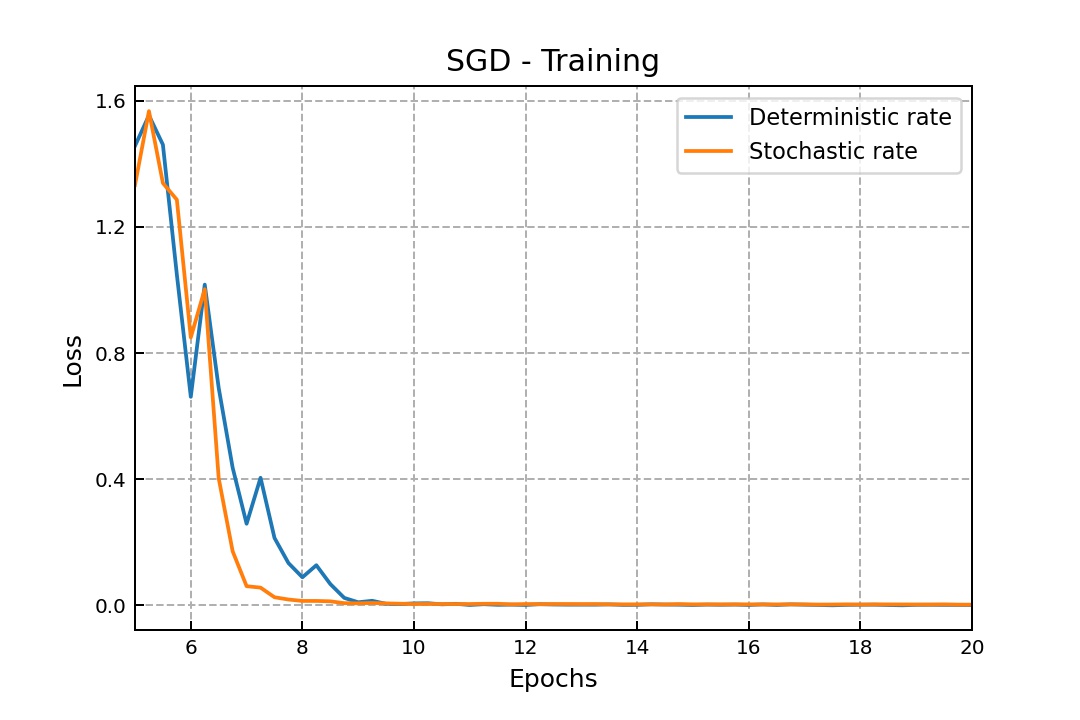}
				\label{fig:}
			\end{subfigure}
			\begin{subfigure}[b]{0.32\textwidth}
				\centering
				\includegraphics[width=\textwidth]{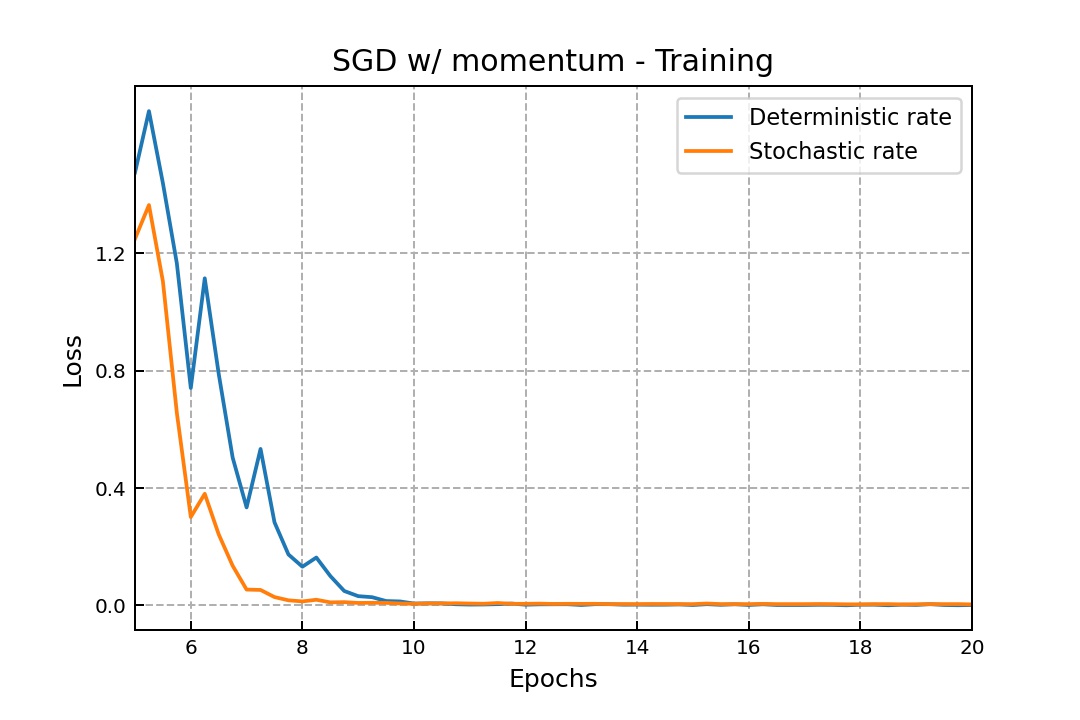}
				\label{fig:}
			\end{subfigure}
			\begin{subfigure}[b]{0.32\textwidth}
				\centering
				\includegraphics[width=\textwidth]{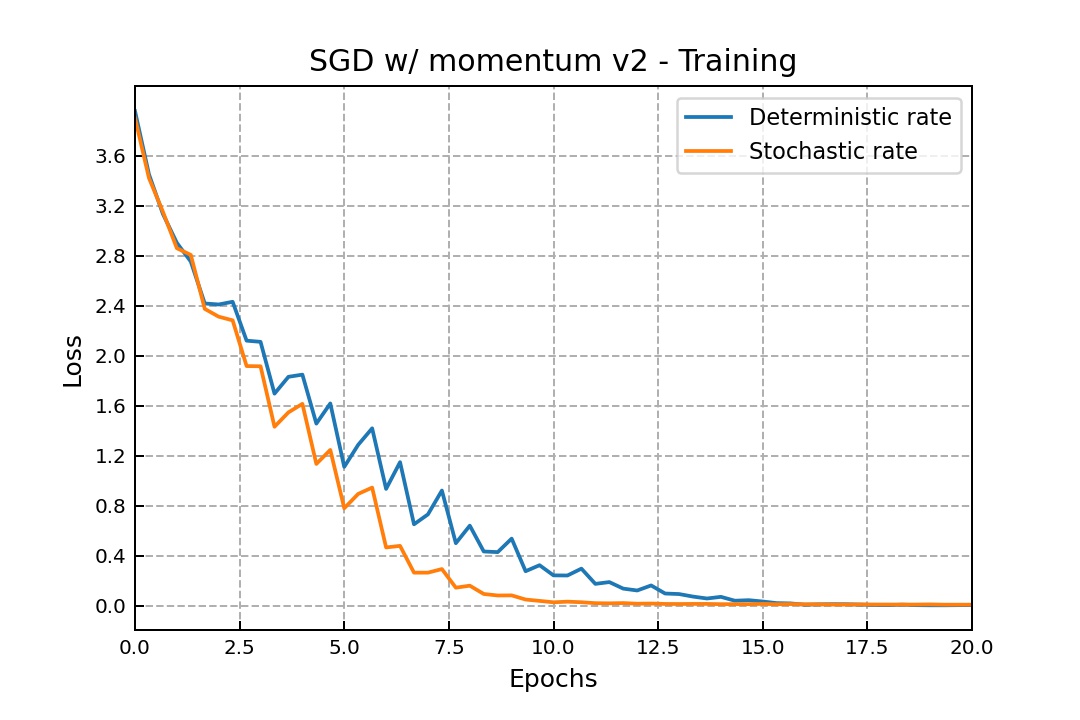}
				\label{fig:}
			\end{subfigure}
			\begin{subfigure}[b]{0.32\textwidth}
				\centering
				\includegraphics[width=\textwidth]{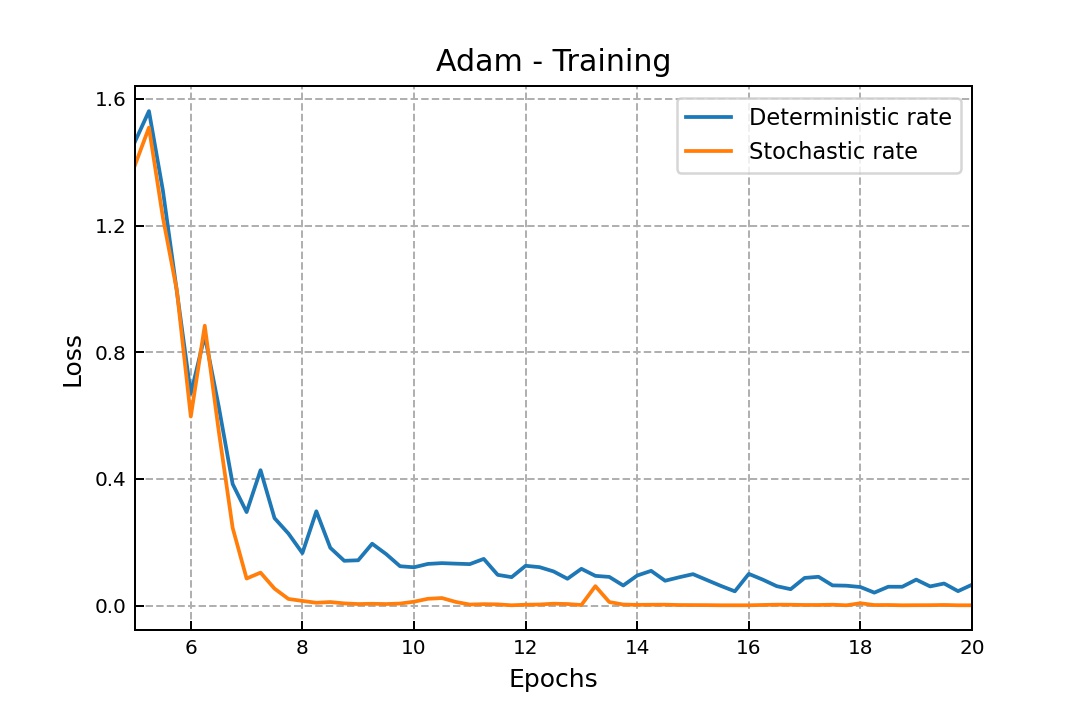}
				\label{fig:}
			\end{subfigure}
			\begin{subfigure}[b]{0.32\textwidth}
				\centering
				\includegraphics[width=\textwidth]{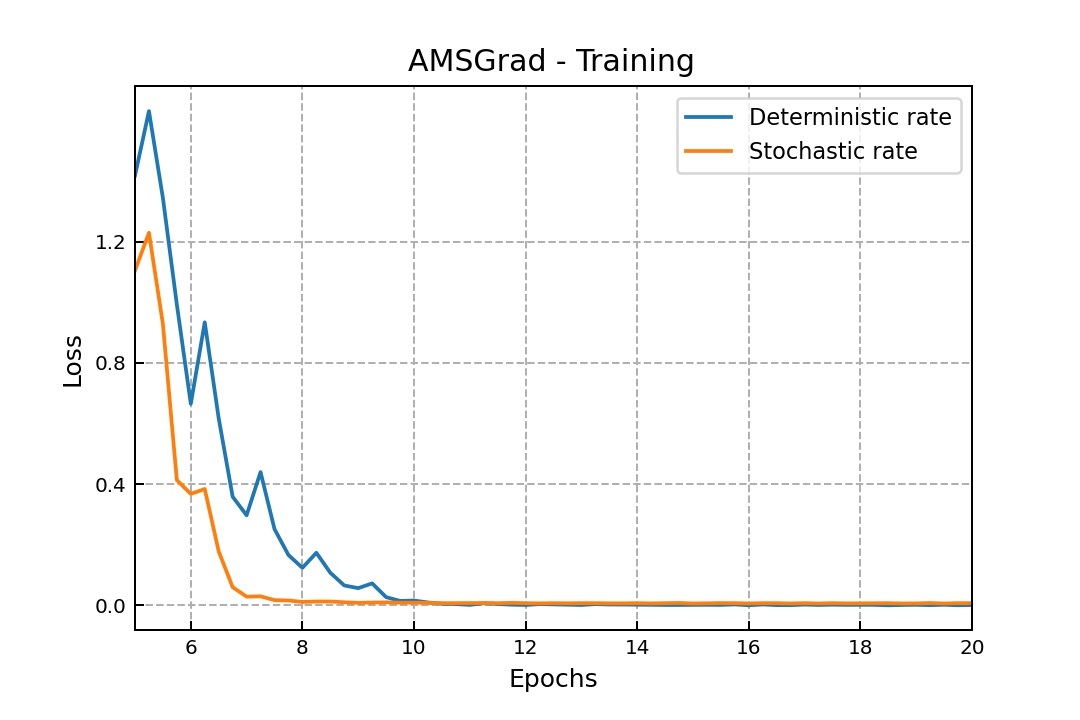}
				\label{fig:}
			\end{subfigure}
			\begin{subfigure}[b]{0.32\textwidth}
				\centering
				\includegraphics[width=\textwidth]{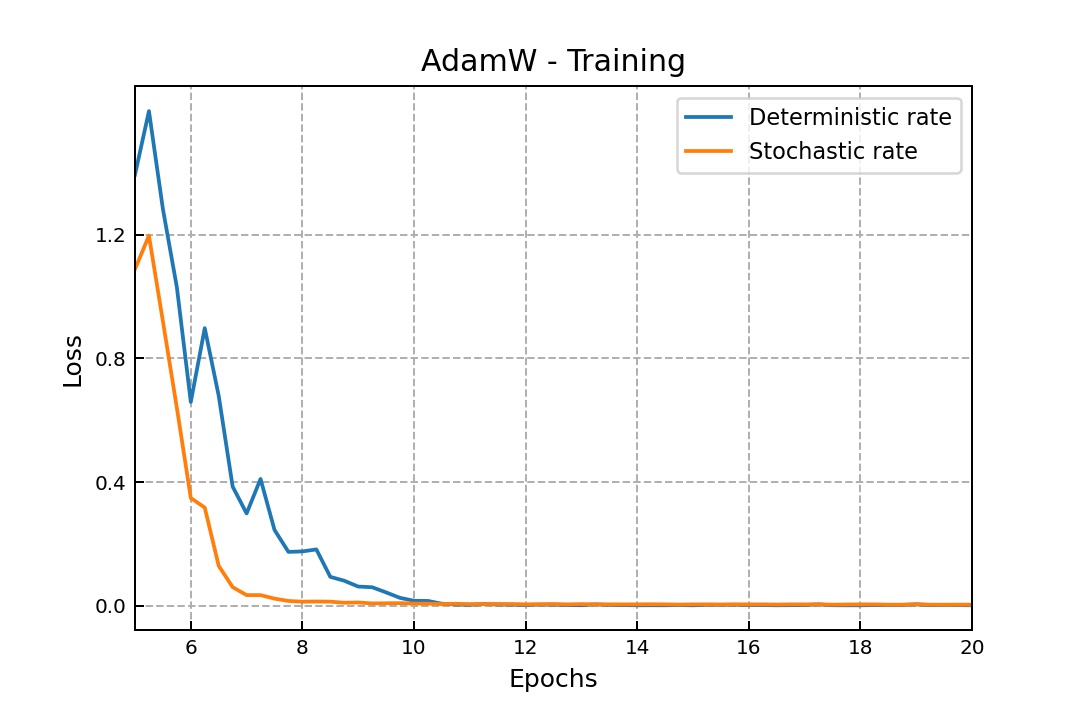}
				\label{fig:}
			\end{subfigure}
			\vspace{-5mm}
			\caption{Plots that illustrate the  $\beta$-RUMSLR (in orange) and the deterministic learning rate (in blue) schemes for $ \beta=100 $, $ c_1=0.9 $ and $ c_2=1.1 $. Both use the same initialization. AMSGrad uses $ c_1=0.85 $, $ c_2=1.15 $. For CIFAR-100, algorithms with  $\beta$-RUMSLR-with-memory yield noticeable minimization performance gains. Where necessary, epochs are omitted for clarity at the start.}
			\label{fig:CIFAR100}
		\end{figure}

		The MNIST dataset was used to validate the online framework theoretical guarantees  for $ \beta-$RUMSLR-with-memory SGD using a convex loss function, in specific logistic regression as aforementioned. Indeed, the results in Fig. \ref{fig:MNIST} verify the convergence of the MNIST $ \beta-$RUMSLR-with-memory scheme in a convex setting as predicted by Theorems \ref{th:theorem_online_R}. Moreover, CIFAR-10 and CIFAR-100, the $ \beta-$RUMSLR-with-memory scheme demonstrated significantly improved convergence than when using a deterministic learning rate  as shown in Figs. \ref{fig:CIFAR10} and \ref{fig:CIFAR100} respectively. As can be seen by the plots, $ \beta-$RUMSLR-with-memory algorithms convergence much faster than their deterministic-learning-rate counterparts.

		\section{Conclusion and Future Work}
		\label{sec:Conclusion and Future Work}
		This work introduced learning-rate schemes that make the learning rate stochastic by multiplicatively equipping it with a function of random variables, which are coined stochastic factors.  Convergence of the SGD algorithm employing stochastic learning rate schemes with memory of past stochastic factors  was theoretically analyzed and compared with known results for the algorithm's deterministic-learning-rate version. The discussion that followed suggested how the hyperparameter values introduced by the two stochastic-learning-rate schemes should be chosen. Empirical results on popular algorithms demonstrated noticeable increase in optimization performance, presenting stochastic-learning-rate schemes as a viable option for enhancing performance. In-depth generalization performance studies, convergence analysis of algorithms besides SGD, and investigating the effect of various stochastic factor distributions and hyperparameters on algorithm performance are some paths for future work.

		\bibliographystyle{abbrvnat}
		\bibliography{refnew_acc2022}
		
		
		\newpage
		\appendix
		\section{APPENDIX}\label{sec:Appendix}
			The appendix contains the proofs for stochastic learning-rate schemes without memory in the stochastic approximation and online learning rate settings. 
			\subsection{Stochastic Optimization Setting for Memoryless Stochastic Learning Rates}
			The theorem that follows provides a bound on the optimality gap of the SGD algorithm when using a memoryless stochastic-learning-rate with a constant step size: 
		\begin{theorem}
			\label{th:theorem_finite_dataset_fixed_step_size_R_old}
			Assume (\ref{eq:L_smooth}-\ref{eq:variance_related_ineq}) are satisfied. Assume a bounded stochasticity factor $ c_1 \le  u_t \le c_2 $.  Furthermore, assume a positive learning rate of:
			\begin{equation} \label{eq:sqrt_lr_R_old}
			a_t=\frac{  c_1}{c_2^2} au_t, 
			\end{equation}
			and:
			\begin{equation}\label{eq:lr_inequality_fixed_step_size_R_sqrt_old}
			a < \frac{c_2}{c_1} \frac{1}{L M_G},
			\end{equation}
			then the iterates of the SGD algorithm in (\ref{eq:SGD}) satisfy:
			\begin{equation} \label{eq:bounds_SGD_constant_step_size_R_sqrt_old}
			\mathop {\lim }\limits_{t \to \infty } {\mathbb{E} }[f({\theta_{t+1}}) - {f_{min}}]  \le 
			\frac{{aLM}}{{2c}}.
			\end{equation}  
		\end{theorem}
		The guarantee in (\ref{eq:bounds_SGD_constant_step_size_R_sqrt_old}) is the same as when SGD uses a deterministic-learning-rate, e.g., as given in \cite{Bottou2018}. The step size is $ \frac{  c_1}{c_2^2} $ times smaller than in the deterministic-learning-rate case. This is pragmatically nonrestrictive for convergence speed as supported by experimental results, since $ c_1 $ and $ c_2 $ are adequately close in practice (see discussion in Section \ref{sec:Discussion on the Stochasticity Factor_old}).

		\begin{proof}
			Taking expectations with respect to $ \xi_t $ in (\ref{eq:L_smooth}), and from (\ref{eq:SGD}) substituting $ \theta_{t+1} $ for $ \theta  $ and $ \theta_t $ for $ \theta ' $, and in (\ref{eq:sqrt_lr_R_old}) letting $ \eta:= \frac{ \sqrt{c_1}}{c_2\sqrt{c_2}} a $ (thus $ a_t=\eta u_t $) yields:
			\begin{equation}
			\begin{aligned}[b] \label{eq:first_expressions_old}
			{\mathbb{E}_{{\xi_t}}}[f({\theta_{t+1}})]  &\le f({\theta_t}) - {a_t}{\nabla ^T}f({\theta_t})\mathbb{E}_{{\xi_t}}[g({\theta_t},{\xi_t})] + a_t^2\frac{L}{2}{\mathbb{E}_{{\xi_t}}}\left[ {{{\left\| {g({\theta_t},{\xi_t})} \right\|}^2}} \right] \\
			& \le   f({\theta_t}) - {a_t}{\left\| {\nabla f({\theta_t})} \right\|^2} + a_t^2\frac{L}{2}{\mathbb{E}_{{\xi_t}}}\left[ {{{\left\| {g({\theta_t},{\xi_t})} \right\|}^2}} \right] \\& \le f({\theta_t}) + {a_t}({a_t}\frac{{{M_G}L}}{2} - 1){\left\| {\nabla f({\theta_t})} \right\|^2} + a_t^2\frac{L}{2}M \\
			&  \le f({\theta_t}) + \frac{1}{2}{a_t}({a_t}{M_G}L - 2){\left\| {\nabla f({\theta_t})} \right\|^2} + a_t^2\frac{L}{2}M  \\
			&  \le f({\theta_t}) + \frac{1}{2}{a_t}({\eta}c_2{M_G}L - 2){\left\| {\nabla f({\theta_t})} \right\|^2} + a_t^2\frac{L}{2}M     \\
			&  \le f({\theta_t}) - \frac{1}{2}{a_t}{\left\| {\nabla f({\theta_t})} \right\|^2} + a_t^2\frac{L}{2}M	\\
			&  \le f({\theta_t}) - c{a_t}(f({\theta_t}) - {f_{min}}) + a_t^2\frac{L}{2}M .
			\end{aligned}	
			\end{equation}
			The second inequality followed from (\ref{eq:unbiased_est_grad}), the third from (\ref{eq:variance_related_ineq}), the fifth from $ a_t=\eta u_t \le \eta c_2 $, the sixth from
			(\ref{eq:lr_inequality_fixed_step_size_R_sqrt_old}) which gives that $ ac_2 M_G L-2  \le -1 $, and the seventh from (\ref{eq:L_smooth_implication}).
			
			Then, using $ \eta c_1\le a_t \le \eta c_2 $ in (\ref{eq:first_expressions_old}),  $ \xi_i $ up to time $ t $ and denoting the result as $ \mathbb{E}_\xi[\cdot] $:
			\begin{equation}
			\begin{aligned}[b] 
			&{\mathbb{E}_\xi }[f({\theta_{t+1}})]  \le {\mathbb{E}_\xi 	}[f({\theta_t})] - c\eta c_1({\mathbb{E}_\xi }[f({\theta_t})] - {f_{min}}) + c_2^2\eta^2\frac{{LM}}{2},
			\end{aligned}	
			\end{equation}
			or equivalently:		
			\begin{equation}
			\begin{aligned}[b] \label{eq:add_and_subtract_R_old}
			{\mathbb{E}_\xi }[f({\theta_{t+1}})]  &\le {\mathbb{E}_\xi }[f({\theta_t})] - c\eta c_1({\mathbb{E}_\xi }[f({\theta_t})] - {f_{min}})\\ & + \eta^2c_2^2\frac{{LM}}{2}  + \eta c_2^2\frac{{LM}}{{2cc_1}} - \eta c_2^2\frac{{LM}}{{2cc_1}}  + {f_{min}} - {f_{min}}.
			\end{aligned}	
			\end{equation}	
			This further results in:
			\begin{equation}
			\begin{aligned}[b]
			{\mathbb{E} }[f({\theta_{t+1}}) - {f_{min}}] -  \eta \frac{{LMc_2^2}}{{2cc_1}} \le (1 - c{{\eta }}c_1)({\mathbb{E} }[f({\theta_t}) - {f_{min}}] - {{\eta}}\frac{{LMc_2^2}}{{2cc_1}}).
			\end{aligned}
			\end{equation}
			By unrolling the inequality:
			\begin{equation}
			\begin{aligned}[b] \label{eq:before_final_fixed_step_size_R_old}
			{\mathbb{E} }[f({\theta_{t+1}}) - {f_{min}}] -  \eta \frac{{LMc_2^2}}{{2cc_1}}  \le (1 - c{{\eta}}c_1)^t({\mathbb{E} }[f({\theta_1}) - {f_{min}}] - {{\eta}}\frac{{LMc_2^2}}{{2cc_1}}).
			\end{aligned}
			\end{equation}
			Furthermore, from (\ref{eq:lr_inequality_fixed_step_size_R_sqrt_old}) it follows that:
			\begin{equation} \label{eq:lr_inequality_eta_fixed_step_size_R_old}
			a < \sqrt{\frac{c_2}{c_1}} \frac{1}{L M_G} \Leftrightarrow \eta = \frac{a \sqrt{c_1}}{c_2\sqrt{c_2}} <\frac{1}{L M_G c_2},
			\end{equation}
			which for $ M_G \ge 1 $, $ c \le L $ yields:
			\begin{equation}
			c \eta c_1 \le c a_t \le c \eta c_2 < \frac{cc_2}{{c_2{M_G}L}} \le  \frac{c}{{L}} \le 1.
			\end{equation}
			Therefore from (\ref{eq:before_final_fixed_step_size_R_old}) it is that:
			\begin{equation}
			\mathop {\lim }\limits_{t \to \infty }	{\mathbb{E} }[f({\theta_{t+1}}) - {f_{min}}] - 
			\eta \frac{{LMc_2^2}}{{2cc_1}} = 0,
			\end{equation}
			or equivalently:
			\begin{equation} \label{eq:last_fixed_step_size_R_old}
			\mathop {\lim }\limits_{t \to \infty } {\mathbb{E} }[f({\theta_{t+1}}) - {f_{min}}]  = 
			\eta \frac{{LMc_2^2}}{{2cc_1}}.
			\end{equation}
			Substituting back $ \frac{ c_1}{c_2^2} a $ for $ \eta $ proves the theorem.
		\end{proof}
		Next, a convergence theorem for diminishing step size is provided.
		\begin{theorem} \label{th:theorem_finite_dataset_decreasing_step_size_R_old}   
			Assume (\ref{eq:L_smooth}-\ref{eq:variance_related_ineq}) are satisfied. Assume a bounded stochasticity factor $ c_1 \le u_t \le c_2 $.
			Furthermore, assume a positive learning rate of $ a_t=\eta_t u_t, $
			Then, if $ \eta_t$ is decreasing, $ \sum\nolimits_{t = 1}^\infty  {\eta _t^2}  < \infty  $ and:
			\begin{equation}\label{eq:lr_decreasing_inequality_R_old}
			\eta_1 < \frac{1}{L M_G c_2},
			\end{equation}
			the iterates of the SGD algorithm in (\ref{eq:SGD}) satisfy:
			\begin{equation} \label{eq:lr_decreasing_final_th_R_old}
			\mathop {\lim }\limits_{t \to \infty } {\mathbb{E} }[f({\theta_{t+1}})] = {f_{min}}.
			\end{equation}     
		\end{theorem}
		The guarantee in (\ref{eq:lr_decreasing_final_th_R_old}) provides the same convergence result as in the deterministic-learning-rate case, it holds for arbitrary decreasing step sizes and does not depend on the initial point unlike previous works (cf. \cite{Bottou2018}).
		\begin{proof}
			The proof for diminishing step size is identical to the proof of Theorem \ref{th:theorem_finite_dataset_fixed_step_size_R_old} until equation (\ref{eq:first_expressions_old}), except that the fifth inequality in (\ref{eq:first_expressions_old}) now follows from $ a_t=\eta_t u_t \le \eta_1 c_2 $, and the sixth from (\ref{eq:lr_decreasing_inequality_R_old}) which gives $ \eta_1 c_2 M_G L-2   \le  -1 $ and therefore $\eta_t c_2 M_G L-2   \le  -1 $, since $ \eta_t \le \eta_1 $ for $ t \in \mathbb{N} $. Then, taking expectations with respect to all $ \xi_i $ up to time $ t $:
			\begin{equation}
			\begin{aligned}[b] \label{eq:pre_add_and_subtract_old}
			{\mathbb{E}_\xi }[f({\theta_{t+1}})]  \le {\mathbb{E}_\xi }[f({\theta_t})] - c{a_t}({\mathbb{E}_\xi }[f({\theta_t})] - {f_{min}}) + a_t^2\frac{{LM}}{2},
			\end{aligned}	
			\end{equation}
			or:
			\begin{equation}
			\begin{aligned}[b]
			{\mathbb{E}_\xi }[f({\theta_{t+1}})]  &\le {\mathbb{E}_\xi }[f({\theta_t})] - c{a_t}({\mathbb{E}_\xi }[f({\theta_t})] - {f_{min}})+ a_t^2\frac{{LM}}{2} + a_t\frac{{LM}}{{2c}} - a_t\frac{{LM}}{{2c}}\\ &   + a_{t+1}\frac{{LM}}{{2c}} - a_{t+1}\frac{{LM}}{{2c}} + {f_{min}} - {f_{min}}.
			\end{aligned}	
			\end{equation}	
			Taking expectations with respect to all learning rate random variables $ a_t $, denoting the result $ \mathbb{E}_a[\cdot] $ and also denoting  $ \mathbb{E}_a[\mathbb{E}_\xi[\cdot]] $ as $ \mathbb{E}[\cdot] $, gives:
			\begin{equation}
			\begin{aligned}[b]
			{\mathbb{E} }[f({\theta_{t+1}})]  &\le {\mathbb{E} }[f({\theta_t})] - c\E[a]{a_t}({\mathbb{E} }[f({\theta_t})] - {f_{min}}) +  \E[a]{a_t^2}\frac{{LM}}{2}  + \E[a]{a_t}\frac{{LM}}{{2c}} - \E[a]{a_t}\frac{{LM}}{{2c}}  \\&+ \E[a]{a_{t+1}}\frac{{LM}}{{2c}} - \E[a]{a_{t+1}}\frac{{LM}}{{2c}}  + {f_{min}} - {f_{min}}.
			\end{aligned}
			\end{equation}
			Then using $ \mathbb{V}_{a}[a_t]=\E[a]{a^2_t}-\E[a]{a_t}^2 $ to replace $ \E[a]{a^2_t} $, and then rearranging terms:
			\begin{equation}
			\begin{aligned}[b]
			{\mathbb{E} }[f({\theta_{t+1}}) - {f_{min}}] - \E[a]{a_{t+1}}\frac{{LM}}{{2c}} &\le (1 - c\E[a]{{a_t}})({\mathbb{E} }[f({\theta_t}) - {f_{min}}] - \E[a]{{a_t}}\frac{{LM}}{{2c}})\\& + (\E[a]{a_t} - \E[a]{a_{t+1}})\frac{{LM}}{{2c}} + \mathbb{V}_a[{a_t}]\frac{{LM}}{2}.
			\end{aligned}
			\end{equation}
			By unrolling the sequence it is:
			\begin{equation}
			\begin{aligned}[b] \label{eq:before_partition_d_old}
			\mathbb{E}[f({\theta_{t+1}}) - {f_{min}}] - {\mathbb{E}_a}[a_{t+1}]\frac{{LM}}{{2c}} &\le \prod\limits_{j=1}^t {(1 - c{\mathbb{E}_a}[{a_j}])} ({\mathbb{E} }[f({\theta _1}) - {f_{min}}] - \E[a]{a_1}\frac{{LM}}{{2c}})\\
			&+ \frac{{LM}}{{2c}}\sum\limits_{j = 1}^{t} {(\E[a]{a_j} - \E[a]{a_{j + 1}})} \prod\limits_{n = j+1}^{t} {(1 - c{\mathbb{E}_a}[{a_{n }}])}   \\ 
			&+ \frac{{LM}}{{2}}\sum\limits_{j = 1}^{t} {\mathbb{V}_{a}[{a_j}]} \prod\limits_{n = j+1}^{t} {(1 - c{\mathbb{E}_a}[{a_{n}}])}.
			\end{aligned}
			\end{equation}	
			Furthermore, from (\ref{eq:lr_decreasing_inequality_R_old}), $ M_G \ge 1 $, and $ c \le L $:
			\begin{equation}
			c a_t \le c a_1  \le c \eta_1 c_2 < \frac{cc_2}{{c_2{M_G}L}} \le  \frac{c}{{L}} \le 1,
			\end{equation}
			it is $ 0 \le 1 - c \mathbb{E}[a_t] < 1 $, since $ c \mathbb{E}[a_t] $ is positive. Moreover, $	c\eta_1c_1 \le  c \eta_1 c_2 <1$. Let $ \delta := 1-c\eta_1 c_1 $.
			Then, it is shown that the last two terms in (\ref{eq:before_partition_d_old}) converge. This is since:
			\begin{align}
			&\sum\limits_{j = 1}^{t} {(\E[a]{a_j} - \E[a]{a_{j + 1}})} \prod\limits_{n = j+1}^{t} {(1- c{\mathbb{E}_a}[{a_{n }}])}  \le { (\eta_1c_2 - \eta_1c_1) }\sum\limits_{j = 1}^t  \prod\limits_{n = j+1}^{t} {(1 - c\eta_1c_1)}  \nonumber \\  &= \eta_1(c_2 - c_1)\sum\limits_{j = 1}^t (1 - c\eta_1c_1)^{t-j-1} \nonumber
			= \eta_1(c_2 - c_1)  \frac{{1 - {\delta ^t}}}{{(1 - \delta )\delta }}  \nonumber,
			\end{align}
			which converges to $ \frac{\eta_1(c_2-c_1)}{(1-\delta)\delta} $ as $ t \to \infty $.
			For the other term:
			\begin{equation}\label{eq:before_partition_old}
			\begin{aligned}[b]
			&\sum\limits_{j = 1}^{t}{\mathbb{V}_{a}[{a_j}]} \prod\limits_{n = j+1}^{t} {(1 - c{\mathbb{E}_a}[{a_{n}}])} \\ & \le \sum\limits_{j = 1}^t \frac{{{\eta_j^2}}}{{4}} {{{(\sqrt[j]{{{c_2}}} - \sqrt[j]{{{c_1}}})}^2}} \prod\limits_{n = j+1}^{t} {(1 - c{\mathbb{E}_a}[{a_{n}}])}  
			\\ & \le \frac{{{\eta_1^2 {{{({c_2} - {c_1})}^2}} }}}{{4}}\sum\limits_{j = 1}^t  \prod\limits_{n = j+1}^{t} {(1 - c\eta_1c_1 )}
			\\ & \le \frac{{{\eta_1^2 {{{({c_2} - {c_1})}^2}} }}}{{4}}\sum\limits_{j = 1}^t \delta^{t-j-1}
			\\ & \le \frac{{{\eta_1^2{{({c_2} - {c_1})}^2}}}}{{4}} \frac{1-\delta^t}{(1-\delta)\delta},
			\end{aligned}
			\end{equation}
			which converges to $  \frac{{{\eta_1^2}{{({c_2} - {c_1})}^2}}}{{4 (1-\delta)\delta }} $ as $ t \to \infty $, and where the second inequality followed from the variance-bound of bounded random variables.
			Then let the algorithm run for long enough, so that the last two terms in (\ref{eq:before_partition_d_old}) are close enough to their limits so that for their tails after some timestep $ K $ the following holds:
			\begin{equation}
			\begin{aligned}[b] \label{eq:aprox_eq_to_zero_d_old}
			\sum\limits_{j = K+1}^{t} {(\E[a]{a_j} - \E[a]{a_{j + 1}})} \prod\limits_{n = j+1}^{t} {(1 - c{\mathbb{E}_a}[{a_{n }}])} & \approx 0 \\    
			\sum\limits_{j = K+1}^{t}{\mathbb{V}_{a}[{a_j}]} \prod\limits_{n = j+1}^{t} {(1 - c{\mathbb{E}_a}[{a_{n}}])} & \approx 0.
			\end{aligned}
			\end{equation}
			Then, partitioning the last two terms in (\ref{eq:before_partition_d_old}) at $ K $, and using (\ref{eq:aprox_eq_to_zero_d_old}), it is that:
			\begin{equation}
			\begin{aligned}[b] \label{eq:partition_d_old}
			\mathbb{E}[f({\theta_{t+1}}) - {f_{min}}] - {\mathbb{E}_a}[a_{t+1}]\frac{{LM}}{{2c}} &\le \prod\limits_{j=1}^t {(1 - c{\mathbb{E}_a}[{a_j}])} ({\mathbb{E} }[f({\theta _1}) - {f_{min}}] - \E[a]{a_1}\frac{{LM}}{{2c}}) \\
			& + \frac{{LM}}{{2c}}\sum\limits_{j = 1}^{K} {(\E[a]{a_j} - \E[a]{a_{j + 1}})} \prod\limits_{n = j+1}^{t} {(1 - c{\mathbb{E}_a}[{a_{n }}])}   \\ 
			& + \frac{{LM}}{{2}}\sum\limits_{j = 1}^{K} {\mathbb{V}_{a}[{a_j}]} \prod\limits_{n = j+1}^{t} {(1 - c{\mathbb{E}_a}[{a_{n}}])}.
			\end{aligned}
			\end{equation}
			Therefore, from (\ref{eq:partition_d_old}), and from $ \E[a]{a_{t+1}}=\eta_{t+1}\E{u_{t+1}} $ it follows that:
			\begin{equation} \label{eq:decreasing_step_size_old}
			\mathop {\lim }\limits_{t \to \infty } \left(	{\mathbb{E} }[f({\theta_{t+1}}) - {f_{min}}] - 
			\eta_{t+1}\E{u_{t+1}}\frac{{LM}}{{2c}} \right) = 0,
			\end{equation}
			which means:
			\begin{equation}
			\mathop {\lim }\limits_{t \to \infty } {\mathbb{E} }[f({\theta_{t+1}}) - {f_{min}}]  =0,
			\end{equation}
			since $  \mathop {\lim }\limits_{t \to \infty }\eta_{t+1}\E{u_{t+1}}\frac{{LM}}{{2c}} =0$ from $  c_1 \le \E{u_{t+1}} \le c_2 $ and $ \eta_{t+1} \to 0 $ as $ t \to \infty $.
		\end{proof}

		Concluding this section, from Theorems
		\ref{th:theorem_finite_dataset_fixed_step_size_R_old}-\ref{th:theorem_finite_dataset_decreasing_step_size_R_old}, it is thus observed that for the stochastic optimization setting the optimality gap bounds for SGD using a stochastic learning rate do not differ from the ones when using a deterministic learning rate. 
		\subsection{Online Learning Setting for Memoryless Stochastic Learning Rates}
		\label{sec:Online Learning Setting_old}
		Now, an alternative analysis for SGD in the case which holds for both MSLR and RMSLR is given for the memoryless case by considering the online learning setting, similarly in Section \ref{sec:Online Learning Setting}.	It will be shown that the memoryless-stochastic-learning-rate SGD has vanishing average regret for diminishing step size $ \eta_t = \frac{a}{\sqrt{t}} $.  The first theorem for the online learning setting follows:
		\begin{theorem}  
			\label{th:theorem_online_R_old}
			Assume that $ X $ is convex, closed, and has bounded diameter $ D $, and that the gradient of $ f_t $ is bounded, i.e. $ \left\| {\nabla {f_t}(\theta )} \right\|_\infty \le G $. Furthermore, assume a positive learning rate of
			$ 	a_t=\frac{a}{\sqrt{t}} u_t, \,\,\, c_1 \le  u_t \le c_2  $
			with $ a = \frac{D}{{G \sqrt {{c_2}{c_1}} }} $.	Then 
			for all $ T\ge 1 $ it is that:
			\begin{equation} \label{eq:online_mslrsgd_beta_old}
			{R_T} \le \sqrt {\frac{{{c_2}}}{{{c_1}}}} 3DG\sqrt T ,
			\end{equation}
		\end{theorem}
		Therefore, $ R_T/T \to 0 $ for $ T \to \infty $.
		\begin{proof}
			Firstly, it is noted that:
			\begin{equation}
			\begin{aligned}[b]
			{\left\| {{\theta _{t + 1}} - {\theta ^*}} \right\|^2} &= {\left\| {{\Pi _X}({\theta _t} - {a_t}{g_t}) - {\theta ^*}} \right\|^2}\\
			&\le {\left\| {{\theta _t} - {a_t}{g_t} - {\theta ^*}} \right\|^2} \\
			&= {\left\| {{\theta _t} - {\theta ^*}} \right\|^2} - 2{a_t}g_t^T({\theta _t} - {\theta ^*}) + a_t^2{\left\| {{g_t}} \right\|^2},
			\end{aligned}
			\end{equation}
			where the first inequality followed from (\ref{eq:projection_property}) replacing $\theta $ by $ \theta^* $ and $ y $ by $ {\theta _t} - {a_t}{g_t} $. Rearranging with respect to $  g_t^T({\theta _t} - {\theta ^*}) $ and using (\ref{eq:f_convex}), it is:
			\begin{equation}
			\begin{aligned}[b]
			{f_t}({\theta _t}) - {f_t}({\theta ^*}) &\le g_t^T({\theta _t} - {\theta ^*})\\
			&\le \frac{{{{\left\| {{\theta _t} - {\theta ^*}} \right\|}^2} - {{\left\| {{\theta _{t + 1}} - {\theta ^*}} \right\|}^2}}}{{2{a_t}}} + \frac{{{a_t}}}{2}{\left\| {{g_t}} \right\|^2}.
			\end{aligned}
			\end{equation}
			Then, summing over all steps:
			\begin{equation} \label{eq:summing_all_old}
			\begin{aligned}[b]
			\sum\limits_{t = 1}^T {[{f_t}(\theta ) - {f_t}({\theta ^*})]}  &\le \sum\limits_{t = 1}^T {\frac{1}{{2{a_t}}}\left( {{{\left\| {{\theta _t} - {\theta ^*}} \right\|}^2} - {{\left\| {{\theta _{t + 1}} - {\theta ^*}} \right\|}^2}} \right)}  + \frac{1}{2}\sum\limits_{t = 1}^T {{a_t}{{\left\| {{g_t}} \right\|}^2}} \\
			&= \frac{1}{{2{a_1}}}{\left\| {{\theta _1} - {\theta ^*}} \right\|^2} + \sum\limits_{t = 2}^T {{{\left\| {{\theta _t} - {\theta ^*}} \right\|}^2}\left( {\frac{1}{{2{a_t}}} - \frac{1}{{2{a_{t - 1}}}}} \right)} + \frac{1}{2}\sum\limits_{t = 1}^T {{a_t}{{\left\| {{g_t}} \right\|}^2}} \\
			& \le \frac{{{D^2}}}{{2{a_1}}} + {D^2}\sum\limits_{t = 2}^T {\left( {\frac{1}{{2{a_t}}} - \frac{1}{{2{a_{t - 1}}}}} \right)}  + \frac{1}{2}\sum\limits_{t = 1}^T {{a_t}{{\left\| {{g_t}} \right\|}^2}} \\
			& \le \frac{{{D^2}}}{{2{a_T}}} + \frac{1}{2}\sum\limits_{t = 1}^T {{a_t}{{\left\| {{g_t}} \right\|}^2}} ,
			\end{aligned}
			\end{equation}
			where the third line followed by using the diameter $ D $ of the convex set $ X $, and the fourth line by using the telescoping sum of the middle term.
			Using $ a_t=\frac{a}{ \sqrt{t} } u_t $ and that $ c_1 \le u_t \le c_2 $ it is:
			\begin{equation}
			\begin{aligned}[b]
			\sum\limits_{t = 1}^T {[{f_t}(\theta ) - {f_t}({\theta ^*})]}  
			&\le \frac{{{D^2}\sqrt T }}{{2a{c_1}}} + \frac{{{ac_2}}}{2}\sum\limits_{t = 1}^T {\frac{{{{\left\| {{g_t}} \right\|}^2}}}{{\sqrt t }}} \\
			&\le \frac{{{D^2}\sqrt T }}{{2a{c_1}}}   + \frac{{{a{G^2}c_2}}}{2}\sum\limits_{t = 1}^T {\frac{1}{{\sqrt t }}}   \le \frac{{\sqrt T }}{{2{c_1}a}}{D^2} + \frac{{a{G^2}{c_2}}}{2}2\sqrt T ,
			\end{aligned}
			\end{equation}
			where the second inequality followed by using the bound on the gradients, and third inequality by using the fact that $\sum\nolimits_{t = 1}^T {\frac{1}{{\sqrt t }}}  \le 2\sqrt T - 1 \le 2\sqrt T $. Then, minimizing the right hand side of the third inequality yields $ a ={D}/({{G\sqrt {{c_2}{c_1}} }})$ for which:
			\begin{equation}
			{R_T} \le \sqrt {\frac{{{c_2}}}{{{c_1}}}} 3DG\sqrt T ,
			\end{equation}
			thus concluding the proof.
		\end{proof}
		\vspace{-1.5mm}
		The next theorem connects SGD with memoryless stochastic-learning-rate schemes with the data dimension $ d $.
		\begin{theorem} 
			\label{th:theorem__online_RE_fixed_step_size_ddepend_old}
			Assume that $ X $ is convex, closed, and has bounded diameter D, and that the gradient of $ f_t $ is bounded, i.e. $ \left\| {\nabla {f_t}(\theta )} \right\|_\infty \le G $. Furthermore, assume a positive learning rate of $ 	a_t=\frac{a}{ \sqrt{t} } u_t, \,\,\, c_1\le u_t \le c_2 $ 
			where $ \left\| {{g_{1:t}}} \right\| $ denotes $\sqrt {\sum\nolimits_{t = 1}^t {{{\left\| {{g_t}} \right\|}^2}} } $. Then, for: 
			\begin{equation} \label{eq:lr_base_ddepend_old}
			a = \frac{1}{{\sqrt {{c_1}{c_2}} }}\frac{D}{{\sqrt{2}\left\| {{g_{1:t}}} \right\|}} 
			\end{equation}
			it is that, for all $ T\ge 1 $ :
			\begin{equation}\label{eq:online_bmslrsgd_old}
			{R_T} \le \sqrt {\frac{{{c_2}}}{{{c_1}}}} \sqrt{2} DG\sqrt {dT} .
			\end{equation}
		\end{theorem}
		\begin{proof}
			The proof is identical to the one of Theorem \ref{th:theorem_online_R_old} until (\ref{eq:summing_all_old}). Then, using $ a_t=\frac{a}{ \sqrt{t} } u_t $  and (\ref{eq:lr_base_ddepend_old}) in (\ref{eq:summing_all_old}) it is:
			\begin{equation}
			\begin{aligned}[b] \label{eq:expressions_before_result_ddependent_old}
			\sum\limits_{t = 1}^T {[{f_t}(\theta ) - {f_t}({\theta ^*})]} & \le \frac{{{D^2}\sqrt {2{c_1}{c_2}} }}{{2D{c_1}}}\left\| {{g_{1:T}}} \right\| + \frac{{{c_2}D}}{{2\sqrt {2{c_1}{c_2}} }}\sum\limits_{t = 1}^T {\frac{{{{\left\| {{g_t}} \right\|}^2}}}{{\left\| {{g_{1:t}}} \right\|}}} \\
			&\le \sqrt {\frac{ {{c_2}}}{{{c_1}}}} \frac{\sqrt{2} D}{2}\left\| {{g_{1:T}}} \right\|  + \sqrt {\frac{{{c_2}}}{{{c_1}}}} \frac{D}{2\sqrt{2} }\sum\limits_{t = 1}^T {\frac{{{{\left\| {{g_t}} \right\|}^2}}}{{\left\| {{g_{1:t}}} \right\|}}} \\
			&\le \sqrt {\frac{{{c_2}}}{{{c_1}}}} \sqrt{2}D\left\| {{g_{1:T}}} \right\|.
			\end{aligned}
			\end{equation}
			The third inequality followed from using $ \sum\nolimits_{t = 1}^T {\frac{{{{\left\| {{g_t}} \right\|}^2}}}{{\left\| {{g_{1:t}}} \right\|}}}  \le 2\left\| {{g_{1:T}}} \right\| $. This can be proved by induction, since for $ t=1 $ it holds that $ \left\| {{g_1}} \right\| \le 2\left\| {{g_1}} \right\| $. Then, assuming the inequality holds for $ T-1 $, for timestep $ T $ it is:
			\begin{equation}
			\begin{aligned}[b]
			\sum\limits_{t = 1}^T {\frac{{{{\left\| {{g_t}} \right\|}^2}}}{{\left\| {{g_{1:t}}} \right\|}}}  = \sum\limits_{t = 1}^{T - 1} {\frac{{{{\left\| {{g_t}} \right\|}^2}}}{{\left\| {{g_{1:t}}} \right\|}} + \frac{{{{\left\| {{g_t}} \right\|}^2}}}{{\left\| {{g_{1:T}}} \right\|}}} & \le 2\left\| {{g_{1:T - 1}}} \right\| + \frac{{{{\left\| {{g_t}} \right\|}^2}}}{{\left\| {{g_{1:T}}} \right\|}} \\ &  \le 2\sqrt {{{\left\| {{g_{1:T}}} \right\|}^2} - g_T^2}  + \frac{{{{\left\| {{g_t}} \right\|}^2}}}{{\left\| {{g_{1:T}}} \right\|}}
			\le 2\left\| {{g_{1:T}}} \right\|.
			\end{aligned}
			\end{equation}
			The last inequality follows from:
			\begin{equation}
			\begin{aligned}[b]
			2\sqrt {{{\left\| {{g_{1:T}}} \right\|}^2} - g_T^2}  & = 2\sqrt {{{\left\| {{g_{1:T}}} \right\|}^2} - \frac{1}{2}\left\| {{g_{1:T}}} \right\|\frac{1}{2}\frac{{g_T^2}}{{\left\| {{g_{1:T}}} \right\|}} + \frac{{g_T^4}}{{4{{\left\| {{g_{1:T}}} \right\|}^2}}} - \frac{{g_T^4}}{{4{{\left\| {{g_{1:T}}} \right\|}^2}}} }  \\
			& = 2\sqrt {{{\left( {{{\left\| {{g_{1:T}}} \right\|}^2} - \frac{{g_T^2}}{{2\left\| {{g_{1:T}}} \right\|}}} \right)}^2} - \frac{{g_T^4}}{{4{{\left\| {{g_{1:T}}} \right\|}^2}}}} \\
			&\le 2\sqrt {{{\left( {{{\left\| {{g_{1:T}}} \right\|}^2} - \frac{{ g_T^2}}{{2\left\| {{g_{1:T}}} \right\|}}} \right)}^2}} 
			\\
			& \le 2{\left\| {{g_{1:T}}} \right\|^2} - \frac{{g_T^2}}{{\left\| {{g_{1:T}}} \right\|}}.
			\end{aligned}
			\end{equation}	
			Then, since  $ g_t:=\nabla f_t(\theta_t) $ for $ \theta_t \in {R^d} $,  the norm $ \left\| {{g_{1:T}}} \right\| $ in (\ref{eq:expressions_before_result_ddependent_old}) is written as $ \sqrt {\sum\nolimits_{t = 1}^T {\sum\nolimits_{j = 1}^d {{{\left| { {{g_{t,j}}} } \right|}^2}} } }  $ where $ g_{t,j} $ is the $ j $th element of the gradient vector at timestep $ t $. Using the bound $ G $ on the gradients it is:
			\begin{equation}
			\sum\limits_{t = 1}^T {[{f_t}(\theta ) - {f_t}({\theta ^*})]}  \le \sqrt {\frac{{{c_2}}}{{{c_1}}}} \sqrt{2}DG\sqrt {dT},
			\end{equation} thus proving the theorem.
		\end{proof}
		In conclusion, from Theorems
		\ref{th:theorem_online_R_old},\ref{th:theorem__online_RE_fixed_step_size_ddepend_old} it is thus observed that for the online optimization setting the regret bounds when using memoryless stochastic learning rates differ by $ \sqrt{\frac{c_2}{c_1}} $ from the bounds of the original SGD algorithm. In specific, the guarantees  in (\ref{eq:online_mslrsgd_beta_old}) and (\ref{eq:online_bmslrsgd_old}) are the same as in the online learning setting for the original SGD, except for a multiplying factor of $ \sqrt {\frac{{{c_2}}}{{{c_1}}}} $  compared to previous works (e.g. \cite{Duchi11}). 
		\subsection{Discussion on the Stochasticity Factor for the Memoryless Learning Rate Scheme}
		\label{sec:Discussion on the Stochasticity Factor_old}
		The bounds in the online learning setting in Theorems \ref{th:theorem_online_R_old},\ref{th:theorem__online_RE_fixed_step_size_ddepend_old} as well as the step size in the stochastic optimization setting in Theorem \ref{th:theorem_finite_dataset_fixed_step_size_R_old} present a trade-off for the choice of $ c_1, c_2 $. On one side they should be kept close enough so that the stochastic-learning-rate SGD step size and online bounds are close to that of the original SGD, since they differ by a factor of $ c_1/c_2^2 $ and $ \sqrt{c_2/c_1} $ respectively. For example $ c_1=0.8,c_2=1.2 $ yield a modest 0.2 increase in the bound of RUMSLR SGD in Theorems \ref{th:theorem_online_R_old},\ref{th:theorem__online_RE_fixed_step_size_ddepend_old} compared to the original SGD, but the further away they are, the larger the bound becomes, and the weaker the guarantees (or smaller the step size in the case of Theorem \ref{th:theorem_finite_dataset_fixed_step_size_R_old}) become. Indeed, it has been observed that choosing low values for $ c_1 $ or large values for $ c_2 $ slows down the algorithm significantly or destabilizes it. On the other side, they should be kept adequately apart so that the bursts are significant enough.

\end{document}